\documentclass{itor}

\usepackage{natbib}%
\usepackage[figuresright]{rotating}

\usepackage{url}
\usepackage{epstopdf}

\usepackage{amsmath}
\usepackage{amssymb}
\usepackage[utf8]{inputenc}
\usepackage{amssymb}
\usepackage{color}
\usepackage{algorithm,algpseudocode}
\usepackage{caption}
\usepackage{array}
\usepackage{float}
\usepackage{comment}
\usepackage{lscape}
\usepackage{booktabs}
\allowdisplaybreaks
\usepackage[titletoc]{appendix}
\usepackage{hyperref}
\usepackage{mathtools}
\usepackage{multirow}
\usepackage{rotating}

\usepackage[english]{babel}
\newtheorem{theorem}{Theorem}

\newenvironment{prf}[1][Proof]{\noindent\textbf{#1.} }{\ \rule{0.5em}{0.5em}}
\newcolumntype{C}[1]{>{\centering\arraybackslash}m{#1}}

\newcommand{\braces}[1]{\left\{#1\right\}}
\newcommand{\reals}{\mathbb{R}}
\newcommand{\expval}[1]{\mathbb{E}\left[ #1 \right]}
\newcommand{\absval}[1]{\lvert #1 \rvert}
\DeclareMathOperator*{\mini}{min.} 
\DeclareMathOperator*{\maxi}{max.} 
\DeclareMathOperator*{\argmin}{arg\,min}
\DeclareMathOperator*{\argmax}{arg\,max}

\jname{International Transactions in Operational Research}
\jvol{XX}
\jyear{20XX}
\doi{xx.xxxx/itor.xxxxx}

\begin{document}

\title{A novel dual-decomposition method for non-convex two-stage stochastic mixed-integer quadratically constrained quadratic problems}

\author[Belyak and Oliveira]{Nikita Belyak\affmark{a} and Fabricio Oliveira\affmark{a,$\ast$}}

\affil{\affmark{a}Department of Mathematics and Systems Analysis, Aalto University, Espoo, Finland}

\email{nikita.belyak@aalto.fi [N. Belyak]; fabricio.oliveira@aalto.fi [F. Oliveira]}

\thanks{\affmark{$\ast$}Author to whom all correspondence should be addressed (e-mail: fabricio.oliveira@aalto.fi).}

\historydate{Received DD MMMM YYYY; received in revised form DD MMMM YYYY; accepted DD MMMM YYYY}

\begin{abstract}
We propose the novel \textit{p}-branch-and-bound method for solving two-stage stochastic programming problems whose deterministic equivalents are represented by non-convex mixed-integer quadratically constrained quadratic programming (MIQCQP) models. The precision of the solution generated by the \textit{p}-branch-and-bound method can be arbitrarily adjusted by altering the value of the precision factor \textit{p}. The proposed method combines two key techniques. The first one, named \textit{p}-Lagrangian decomposition, generates a mixed-integer relaxation of a dual problem with a separable structure for a primal non-convex MIQCQP problem. The second one is a version of the classical dual decomposition approach that is applied to solve the Lagrangian dual problem and ensures that integrality and non-anticipativity conditions are met once the optimal solution is obtained. This paper also presents a comparative analysis of the \textit{p}-branch-and-bound method efficiency considering two alternative solution methods for the dual problems as a subroutine. These are the proximal bundle method and Frank-Wolfe progressive hedging. The latter algorithm relies on the interpolation of linearisation steps similar to those taken in the Frank-Wolfe method as an inner loop in the classic progressive hedging. The \textit{p}-branch-and-bound method's efficiency was tested on randomly generated instances and demonstrated superior performance over commercial solver Gurobi.
\end{abstract}


\maketitle

\newpage

\section{Introduction}
\label{sec: Intro}
Presently, the majority of engineering sectors utilise mathematical optimisation as a modelling framework to represent the behaviour of various processes. Areas such as electrical and process engineering are arguably the most prominent employers of mathematical optimisation techniques to improve their operational performance. For instance, \cite{andiappan_2017} emphasises the efficiency of mathematical optimisation as an approach for designing energy systems. The author highlights the superior performances of mathematical optimisation in terms of comprehensiveness and ability to explicitly determine the topology of energy systems compared to its alternatives, e.g., heuristic and insight-based approaches. \cite{GROSSMANN2019474} discuss the importance of mathematical optimisation in chemical and petrochemical operations, allowing, in some cases, up to 30\% in energy savings. 

Mathematical optimisation approaches closely rely on solving to optimality a mathematical optimisation problem - the set of mathematical relationships representing the real-world problem \citep{kallrath_optimization:_2021}.  \cite{kallrath_optimization:_2021} classifies mathematical optimisation problems into four groups such as linear programming problems, mixed-integer linear programming (MILP) problems, nonlinear programming problems, and mixed-integer nonlinear programming (MINLP) problems. Mixed-integer problems involve decision variables that can have both continuous and discrete domains. The linearity or nonlinearity of the problem refers to the type of constraints and objective function. \cite{lee_mixed_2011} highlight that MINLP problems are particularly challenging due to the difficulties arising from solving over integer variables and non-linear functions. At the same time, both mixed-integer linear programming and nonlinear programming problems are known to be NP-hard \citep{forrest_practical_1974, murtagh1995minos}. Nevertheless, the range of applications of MINLP is noticeably diverse \citep{lee_mixed_2011}. It includes modelling block layout design problems with unequal areas \citep{castillo2005optimization}, structural flow sheet problem \citep{kocis1987relaxation} and finding optimal design of water distribution networks \citep{bragalli2012optimal} to mention only a few relevant applications. 

In this paper, we focus on a subclass of MINLP problems that represent deterministic equivalents for two-stage stochastic mixed-integer programming (2SSMIP) problems. Such problems involve two decision-variable sets that are separated by an intermediate probabilistic event. These two distinct decision-variable sets represent decisions made at different stages, i.e., before and after the intermediate probabilistic event occurred and its outcome is observed. The modelling of probabilistic events for the most part involves consideration of mutually exclusive and exhaustive alternatives (scenarios) and the definition of probabilities associated with them \citep{bush1953stochastic}. Despite their vast applicability, 2SSMIP problems raise serious conceptual and computational challenges \citep{kuccukyavuz2017introduction}. For instance, in \cite{LIAO2019265}, the authors exploited a multi-step mixed-integer nonlinear programming problem to optimise the recovery process for network and load during power system restoration. In \cite{WANG2021143}, a 2SSMIP model has been used to formulate a container slot allocation problem for a liner shipping service. The authors used the sample-average approximation to approximate the expected value function rendering a nonlinear integer programming model. 

Our interest resides in 2SSMIP problems whose deterministic equivalent representations are non-convex mixed-integer quadratically constrained quadratic programming (MIQCQP) models. These models arise in several practically relevant applications, such as the pooling problem, which is a non-convex MIQCQP under the assumption of linearly blending qualities \citep{misener2012global} and as the equivalent, single-level reformulation of some bi-level optimisation problems \citep{ding2014bi, virasjoki2020utility}. The examples of pooling problem applications include the design of water networks \citep{jezowski2010review}, modelling refinery processes \citep{amos1997modelling}, and transportation systems \citep{calvo2004distributed} as some relevant examples. 

The formulation of a general 2SSMIP is 
\begin{equation}
 z^\text{SMIP} = \mini_x \braces{c^\top x + \mathcal{Q}(x) : x \in X} \label{eq:SMIP},
\end{equation}
where the vector $c \in \reals^{n_x}$ is known, $X$ is a mixed-integer linear set consisting of linear constraints and integrality restrictions on some components of $x$. The recourse function $\mathcal{Q} : \reals^{n_x} \to \reals$ is the expected recourse value function
\begin{equation}
 \mathcal{Q}(x) = \expval{\mini_y \braces{ f(y,\xi) : g(x,y,\xi) = 0, y \in Y(\xi)}},
\end{equation}
where, for any realisation of the random variable $\xi$, $f: \reals^{n_y} \to \reals$ is defined as 
$$
f(y, \xi) = q(\xi)^\top y + \sum_{(i,j) \in B_Q} Q(\xi)_{i,j}y_iy_j,
$$ 
$g = [g_1, \dots, g_{ \absval{M}}]^\top$ where $g_m: \reals^{n_x \times n_y} \to \reals, \forall m \in \braces{1,\dots, \absval{M}} = M$, is defined as 
$$
g_m(x,y,\xi) = T(\xi)_mx + W(\xi)_my + \sum_{(i,j) \in B_U} U(\xi)_{m,i,j}y_iy_j - h(\xi)_m,
$$ 
and $B_Q$ ($B_U$) comprise the index pairs $(i,j)$ for which the entry $|Q_{i,j}| > 0$ ($|U_{i,j}| > 0$), implying the presence of the bi-linear terms $y_i y_j$; $Y(\xi)$ is a mixed-integer set containing both linear constraints and integrality requirements on some of the variables $y(\xi)$; and $\expval{\, \cdot \,}$ denotes the expectation of $~\cdot~$ in terms of the random variable $\xi$. As it is standard practice in the stochastic programming literature, we assume that the random variable $\xi$ is represented by a finite set $S$ of realisations $\xi_1, \dots, \xi_{\absval{S}}$, each with associated probability value $\pi_1, \dots, \pi_{\absval{S}}$. In particular, each realisation $\xi_s$ of $\xi$ encodes the realisation observed for each of the random elements $(q(\xi_s), Q(\xi_s))$ and $(T(\xi_s)_m, W(\xi_s)_m, U(\xi_s)_m, h(\xi_s)_m)$, $\forall m \in M$. For the sake of notation compactness, we refer to these collections as $(q^s, Q^s)$ and $(T_m^s, W_m^s, U_m^s, h_m^s)$, $\forall m \in M$, respectively. 

Problem \eqref{eq:SMIP} can be then posed as the deterministic equivalent 
\begin{equation}
\label{eq: DEM}
\begin{aligned}
    z^\text{SMIP} = \mini_{x,y} ~& c^\top x + \sum_{s \in S} \pi^s (q^{s\top} y^s + \sum_{(i,j) \in B_Q} Q^s_{i,j} y^s_i y^s_j)  \\
    \text{subject to:} ~& x \in X \\
    & T_m^s x + W_m^s y^s + \sum_{(i,j) \in B_U} U_{m,i,j}^s y^s_i y^s_j = h_m^s, \ \forall m \in M, \forall s \in S \\
    & y^s \in Y^s, \ \forall s \in S. 
\end{aligned} 
\end{equation}

Due to the challenging nature of the non-convex MIQCQP problems open source and commercial global solvers, such as Gurobi \citep{gurobi}, Couenne \citep{couenne}, or Baron \citep{baron} still present performance issues when addressing large-scale instances. There have been several solution approaches developed for non-convex MIQCQP problems, which can be categorised into three groups. The first one involves approximation of the problem \eqref{eq: DEM} with a continuous or mixed-integer relaxation \citep{cui2013convex, andrade2018rnmdt}. Another group is formed by those employing variants of the branch-and-bound (BnB) method. In particular, typically for non-convex problems, spatial BnB is used, which involves convexification of non-convex terms as a sub-routine \citep{castro2016spatial, ding2014bi, berthold2012extending}. The last group involves methods relying on the introduction of non-anticipativity conditions (NAC) and the decomposition of the problem into more tractable sub-problems.

The block-angular structure of the problem \eqref{eq: DEM} allows for formulating an almost decomposable equivalent problem by making explicit non-anticipativity conditions (NAC) of the first-stage variables $x$. The reformulated deterministic equivalent model (RDEM) with an almost separable structure can be represented as
\begin{equation}
\label{eq: RDEM}
\begin{aligned}
    z^\text{SMIP} = \mini_{x,y} ~&  \sum_{s \in S} \pi^s (c^\top x^s + q^{s\top} y^s + \sum_{(i,j) \in B_Q} Q^s_{i,j} y^s_i y^s_j)  \\ 
    \text{s.t.:} ~&  y^s \in Y^s, \ \forall s \in S \\ 
    & x^s \in X, \ \forall s \in S \\
    & T_m^s x^s + W_m^s y^s +  \sum_{(i,j) \in B_U} U_{m,i,j}^s y^s_i y^s_j = h_m^s, \ \forall m \in M, \ \forall s \in S \\
    & x^s - \overline{x} = 0, \ \forall s \in S, 
\end{aligned} 
\end{equation}
where the constraint $ x^s - \overline{x} = 0, \ \forall s \in S $, enforces non-anticipativity for the first-stage decisions and  $\overline{x}$ is an auxiliary variable used to enforce nonanticipativity for all $x^s, \ \forall s \in S $. The RDEM problem \eqref{eq: RDEM} could be fully decomposed into $\absval{S}$ non-convex MIQCQP problems if one could remove the set of linear constraints $ x^s - \overline{x} = 0, \ \forall s \in S $, that relates variables from distinct sub-problems, a structure commonly known as complicating constraints.  

To tackle the non-convex problem \eqref{eq: RDEM}, \cite{p_lagrangian_Andrade} developed an algorithm named \textit{p}-Lagrangian decomposition. The \textit{p}-Lagrangian decomposition method involves exploiting Lagrangian relaxation for decomposing the primal problem \eqref{eq: RDEM} into $\absval{S}$ independent sub-problems and employing the reformulated normalised multiparametric disaggregation technique (RNMDT) \citep{andrade2018rnmdt} to construct mixed-integer-based relaxations. The necessity to formulate a mixed-integer-based relaxation arises from the non-convex nature of the primal problem \eqref{eq: RDEM}, leading to a lack of monotonicity in the behaviour of its dual bound value. Therefore, relying solely on Lagrangian decomposition becomes insufficient in ensuring a valid lower (dual) bound for the primal problem \eqref{eq: RDEM}. However, applying Lagrangian decomposition to the mixed-integer-based relaxation of the primal problem yields a valid lower (dual) bound. As a subroutine, the \textit{p}-Lagrangian decomposition algorithm employs a dynamic precision-based method, ensuring the tightening of the relaxation bounds as the precision parameter value approaches $-\infty$, and a bundle method approach for updating the dual multipliers. Additionally, the decomposable structure of the Lagrangian dual problem is amenable to parallelisation, which can significantly enhance the computational performance.

As suggested by the numerical results in \cite{p_lagrangian_Andrade}, the \textit{p}-Lagrangian decomposition demonstrated superior performance compared to commercial solver Gurobi \citep{gurobi} when solving non-convex MIQCQP problems with the decomposable structure. Nevertheless, the \textit{p}-Lagrangian decomposition algorithm has an important shortcoming related to the duality gap arising from the mixed-integer nature of the primal problem combined with the imprecision of the RNMDT relaxation.  It is worth highlighting that the convergence of a $p$-Lagrangian relaxation problem requires the precision parameter value to approach $-\infty$. Nevertheless, for most of the practical applications convergence with a predetermined tolerance is sufficient. 

Our primary contribution involves the introduction of a solution method named \textit{p}-branch-and-bound (\textit{p}-BnB) for non-convex MIQCQP problem \eqref{eq: RDEM}. The \textit{p}-BnB algorithm evolves from the advancement of \textit{p}-Lagrangian relaxation, effectively mitigating the duality-gap issue and ensuring convergence to a global optimum. Similar to the performance observed with the \textit{p}-Lagrangian method, the proposed \textit{p}-BnB surpasses the performance of the commercial solver Gurobi \cite{gurobi}. Importantly, our methodology marks the first instance of incorporating \textit{p}-Lagrangian relaxation within the framework of a duality-based branch-and-bound approach. The \textit{p}-BnB method is inspired by the decomposition method for two-stage stochastic integer programs proposed in \cite{caroe1999dual}. The technically challenging synchronisation of \textit{p}-BnB and decomposition method proposed in \cite{caroe1999dual}  relies on the repeatedly solving \textit{p}-Lagrangian relaxation of problem \eqref{eq: RDEM} by means of \textit{p}-BnB and iteratively restricting the feasible region via branch-and-bound framework whenever the solution of \textit{p}-Lagrangian relaxation violates integrality or non-anticipativity conditions. Consequently, \textit{p}-BnB provides the upper bound for problem \eqref{eq: RDEM} that can be made arbitrarily precise against the value of the Lagrangian relaxation bound by decreasing the value of precision factor \textit{p}.

Our subsequent contribution involves the evaluation of the numerical efficiency of \textit{p}-BnB on randomly generated instances considering two different solution methods for the dual problem stemming from \textit{p}-Lagrangian relaxation. The first one is the Frank-Wolfe Progressive Hedging (FWPH) method, originally presented in \cite{FW_PH_2018}. The classic progressive hedging \cite{rockafellar1991scenarios} method is proved to converge to the solution of the deterministic equivalent of a two-stage stochastic programming if such a point exists and the problem is convex. However, there is no guarantee of convergence in case the primary problem is mixed-integer, hence non-convex. In contrast, the FWPH is an enhancement of the classic progressive hedging method with convergence guarantees even in the presence of mixed-integer variables. That is, FWPH guarantees convergence to the optimal dual value of \textit{p}-Lagrangian relaxation. The other solution method for dual problems tested in \textit{p}-BnB is the proximal bundle method \citep{oliveira2014bundle, kim2022scalable}. The proximal bundle method relies on the classic bundle method \citep{lemarechal1974algorithm} but involves regularisation of the dual space search allowing for fewer iterations until convergence \cite{kim2022scalable}. 

The first step of the proposed \textit{p}-BnB method involves the construction of the mixed-integer relaxation of the primal non-convex RDEM problem \eqref{eq: RDEM} by means of employing the RNMDT technique described in Section \ref{sec:RNMDT} to relax all quadratic terms. Next, we apply a Lagrangian duality-based branch-and-bound method reviewed in Section \ref{sec:dual_bnb}. The method involves the composition of branch-and-bound strategies and dual decomposition. Subsequently, each branching tree node is represented by an assembly of decomposable dual sub-problems. To solve the dual sub-problems within the branch-and-bound search, we consider the FWPH method discussed in Section \ref{sec: FWPH} and the proximal bundle method presented in Section \ref{sec: bundle_method}. To the best of our knowledge, this is the first time the efficiency of the FWPH method has been assessed within a Lagrangian duality-based branch-and-bound framework. The proposed method was tested on randomly generated instances, and the results of the numerical experiments are presented in Section \ref{sec:numerical experiments}. Finally, in Section \ref{sec:conclusion}, we provide conclusions and directions for further research.

\section{Technical background}

In what follows, we present the technical elements that form our proposed method. In essence, \textit{p}-BnB is formed by the combination of three main techniques, namely \textit{p}-Lagrangian decomposition, of which RNMDT is a key concept, solution methods for dual Lagrangian problems (FWPH and proximal bundle method), and a branch-and-bound coordination algorithm.   

\subsection{Reformulated normalized multiparametric disaggregation technique (RNMDT)}
\label{sec:RNMDT}

To ensure the validity of the dual (lower) bound value for the non-convex problem \eqref{eq: RDEM} we utilise arbitrarily precise mixed-integer linear relaxation of \eqref{eq: RDEM}, which will then serve as the basis for deriving the corresponding dual problem.  

The normalised multiparametric disaggregation technique (NMDT) is an efficient technique to relax quadratic terms in non-convex MIQCQP problems \citep{castro2016normalized}. NMDT involves the discretisation of the domain of one variable in each bi-linear term. The discretisation procedure in NMDT closely resembles the piecewise McCormick envelopes approach \citep{BERGAMINI20051914, mccormic_1976} as it involves splitting the variable domain into a number of uniform partitions. The accuracy of the NMDT approximation of the primary non-convex MIQCIP problem is directly related to the size of these partitions and can be made arbitrarily precise. However, improving the accuracy of NMDT approximation follows the increase in the number of continuous and binary variables. Therefore, there is a trade-off between the accuracy of the NMDT approximation and its computational tractability. 

The enhancement of NMDT has been proposed by \cite{andrade2018rnmdt} and this new version was named reformulated NMDT (RNMDT). The authors suggested reducing the numerical base used for the discretised domain from 10 to 2 (or binary). Hence, the smallest interval used for discretising the domain would be 0.5 when $p = -1$ instead of 0.1 as in the case of base 10. This modification allowed for a significant reduction in the number of auxiliary binary variables required in NMDT relaxation. Additionally, the authors performed a series of reformulations leading to the elimination of a number of redundant constraints and variables. Therefore, RNMDT relaxation became more easily tractable compared to NMDT. For further details on the reformulation of NMDT relaxation please refer to \cite{andrade2018rnmdt}. The RNMDT relaxation of the primal RDEM problem can be constructed by employing RNMDT to discretise the second-stage variables $y^s_j$ in the primal RDEM. Following the notation in \cite{andrade2018rnmdt}, formally, for a fixed value of the precision factor \textit{p}, the mixed-integer relaxation RNMDT$_p$ can be stated as

\begin{equation}
\label{RNMDT_p}
\begin{aligned}  
    z^{\text{RNMDT}} = &\mini_{x,y,w}  \sum_{s \in S} \pi^s (c^\top x^s + q^{s\top} y^s + \sum _{(i,j) \in B_Q} Q^s_{i,j} w^s_{i,j}) \\ 
    \text{s.t.: } ~&  y^s \in Y^s, \ \forall s \in S, \\
    & x^s \in X, \ \forall s \in S \\
    & x^s - \overline{x} = 0, \ \forall s \in S \\
    & T_m^s x^s + W_m^s y^s + \sum_{(i,j) \in B_U}U_{m,i,j}^s w^s_{i,j} = h_m^s, \ \forall m \in M, \forall s \in S \\
    & y^s_{j} = (N^{U,s}_j - N^{L,s}_j)\bigg( \sum_{l \in P} 2^l z^s_{j, l} + \Delta y^s_j \bigg) + N^{L,s}_j, \ \forall s \in S, \ \forall j \in \{j\mid (i,j) \in B_Q \cup B_U \} \\
    & 0 \le \Delta y^s_{j} \le 2^{p}, \ \forall s \in S, \ \forall j \in \{j\mid (i,j) \in B_Q \cup B_U \} \\ 
    & w^s_{i, j} = (N^{U,s}_j - N^{L,s}_j)\bigg(\sum_{l \in P} 2^l \hat{y}^s_{i, j , l} + \Delta w^s_{i, j} \bigg) + y^s_i N^{L,s}_j, \ \forall s \in S, \ \forall (i, j) \in B_Q \cup B_U \\
    \ \ \ \ \ \ \ \ \ \ \ \ & 2^{p}(y^s_{i}-N^{U,s}_i) + N^{U,s}_i \Delta y^s_{j} \le \Delta w^s_{i, j} \le 2^{p}(y^s_{i}-N^{L,s}_i) + N^{L,s}_i\Delta y^s_{j}, \ \forall s \in S, \ \forall (i, j) \in B_Q \cup B_U  \\
   & N^{L,s}_i \Delta y^s_{j} \le \Delta w^s_{i, j} \le N^{U,s}_i \Delta y^s_{j}, \ \forall s \in S, \ \forall (i, j) \in B_Q \cup B_U  \\
   & N^{L,s}_i z^s_{j, l} \le \hat{y}^s_{i, j, l} \le N^{U,s}_i z^s_{j, l},  \ \forall s \in S, \ \forall (i,j) \in B_Q \cup B_U , l \in P \\
    & N^{L,s}_i (1 - z^s_{j, l}) \le y^s_{i} - \hat{y}^s_{i, j, l} \le N^{U,s}_i (1 - z^s_{j, l}), \ \forall s \in S, (i,j) \in B_Q \cup B_U , l \in P\\ 
     & z^s_{j, l} \in \{0, 1\}, \ \forall s \in S, \ \forall j \in \{j\mid (i,j) \in B_Q \cup B_U \}, \ \forall l \in P,  
\end{aligned} 
\end{equation}
where $P = \{p, \dots, -1\}$ and $w^s_{i,j}$ represent the product $y^s_i y^s_j$.\cite{andrade2018rnmdt} have demonstrated that as the precision parameter $p$ approaches $-\infty$ the corresponding RNMDT$_p$ relaxation becomes increasingly tighter.

\subsection{\textit{p}-Lagrangian relaxation}

Let us consider the RNMDT$_p$ problem \eqref{RNMDT_p} defined in Section \ref{sec:RNMDT}, where the precision factor \textit{p} is fixed to some negative integer value. The \textit{p}-Lagrangian decomposition of  RNMDT$_p$ can be obtained by applying Lagrangian relaxation to relax the NAC 
$$ 
x^s - \overline{x} = 0, \ \forall s \in S.
$$
Let $\lambda = (\lambda^1, \dots, \lambda^{\mid S \mid}) \in \mathbb{R}^{n_x \times \absval{S}} $ be the vector of dual multipliers associated with the relaxed NAC. By setting $\mu^s = \frac{1}{\pi^s} \lambda^s$, $\forall s \in S$, the \textit{p}-Lagrangian dual function can be defined as  
\begin{align}
L(\mu) = \left\{ \begin{aligned}
 \mini_{x,\overline{x},y,w} \sum_{s \in \mathcal{S}}\  \pi^s \big(c^\top x^s + q^{s\top} y^s + \sum _{(i,j) \in B_Q} Q^s_{i,j} w^s_{i,j} + {\mu^s}^\top (x^s - \overline{x}))  \\ : (x^s,y^s, \Gamma^s) \in G^s, \forall s \in S 
\end{aligned} \right\},
\label{p_LR}
\end{align}
where $\Gamma^s = \{ w^s, \Delta y^s, \Delta w^s, \hat{y}^s, z^s\}$ and $G^s$ is defined by the following set of constraints
\begin{equation*}
G^s =\begin{cases} 
     x^s \in X  \\
     y^s \in Y^s \\
     T_m^s x + W_m^s y^s + \sum_{(i,j) \in B_U}U_{m,i,j}^s w^s_{i,j} = h_m^s, \ \forall m \in M \\
     y^s_{j} = (N^{U,s}_j - N^{L,s}_j)\bigg( \sum_{l \in P} 2^l z^s_{j, l} + \Delta y^s_j \bigg) + N^{L,s}_j, \ \forall j \in \{j\mid (i,j) \in B_Q \cup B_U \} \\ 
     0 \le \Delta y^s_{j} \le 2^{p}, \ \forall j \in \{j\mid (i,j) \in B_Q \cup B_U \} \\  
     w^s_{i, j} = (N^{U,s}_j - N^{L,s}_j)\bigg(\sum_{l \in P} 2^l \hat{y}^s_{i, j , l} + \Delta w^s_{i, j} \bigg) + y^s_i N^{L,s}_j, \ \forall (i, j) \in B_Q \cup B_U  \\ 
     2^{p}(y^s_{i}-N^{U,s}_i) + N^{U,s}_i \Delta y^s_{j} \le \Delta w^s_{i, j} \le 2^{p}(y^s_{i}-N^{L,s}_i) + N^{L,s}_i\Delta y^s_{j}, \  \forall (i, j) \in B_Q \cup B_U \\  
    N^{L,s}_i \Delta y^s_{j} \le \Delta w^s_{i, j} \le N^{U,s}_i \Delta y^s_{j}, \ \forall (i, j) \in B_Q \cup B_U  \\ 
    N^{L,s}_i z^s_{j, l} \le \hat{y}^s_{i, j, l} \le N^{U,s}_i z^s_{j, l}, \ \forall (i,j) \in B_Q \cup B_U , \ \forall l \in P \\ 
     N^{L,s}_i (1 - z^s_{j, l}) \le y^s_{i} - \hat{y}^s_{i, j, l} \le N^{U,s}_i (1 - y^s_{j, l}), \ \forall (i,j) \in B_Q \cup B_U , \ \forall l \in P\\ 
      z^s_{j, l} \in \{0, 1\}, \ \forall j \in \{j\mid (i,j) \in B_Q \cup B_U \}, \ \forall l \in P.
\end{cases} 
\end{equation*}

The variable $\overline{x}$ in \eqref{p_LR} is unconstrained. Therefore, in order for the \textit{p}-Lagrangian dual function $L(\mu)$ to be bounded, we must impose the dual feasibility condition $\sum_{s \in S} \pi^s\mu^s =0$. With this assumption in mind, the \textit{p}-Lagrangian dual function \eqref{p_LR} can be explicitly decomposed for each scenario $s \in S$
\begin{align}
L(\mu) = \sum_{s \in \mathcal{S}}\  \pi^s L^s(\mu^s), \label{eq: p_LR_decomposed}
\end{align}
where
\begin{align}
 L^s(\mu^s) = \left\{ \begin{aligned} \mini_{x,y,w} {(c^s + \mu^s)}^{\top} x^s + q^{s^{\top}} y^s + \sum_{(i,j) \in B_Q} Q^s_{i,j} w^s_{i,j} \\ : (x^s,y^s, \Gamma^s) \in G^s, \forall s \in S
 \end{aligned} \right\}.
\label{eq: p_LR_decomposed_subproblem} 
\end{align}

For any fixed value of $\mu = (\mu_1, \dots \mu_{\mid S \mid})$, the \textit{p}-Lagrangian dual function \eqref{eq: p_LR_decomposed} provides a lower bound for the primal RNMDT$_p$ problem \eqref{RNMDT_p} \citep{p_lagrangian_Andrade}. Our objective is to find the tightest (i.e., the uppermost) lower bound. Therefore, the dual bound can be obtained by solving the \textit{p}-Lagrangian dual problem
\begin{align}
    z_{LD} = \max_{\mu} \left\{L(\mu): \sum_{s \in S} \pi^s \mu^s = 0 \right\}.
    \label{eq: dual_problem}
\end{align}

\subsection{Solution method for \textit{p}-Lagrangian dual function}

In this section, we present adaptations of proximal bundle method \citep{kim2022scalable} and FWPH \citep{FW_PH_2018} for solving the  \textit{p}-Lagrangian dual problem \eqref{eq: dual_problem}. One should bear in mind that alternative nonsmooth (convex) optimisation algorithms can be potentially applied to solve \textit{p}-Lagrangian dual function $L(\mu)$. The choice for proximal bundle method and FWPH was motivated by the literature on dual Lagrangian-based methods, including their reported efficiency (see, for example, \cite{FW_PH_2018, palani2019frank, bashiri2021two, feltenmark2000dual}), and our own experience with preliminary experiments involving both and other simpler nonsmooth optimisation methods (i.e., subgradient and cutting planes methods).

\subsubsection{Proximal bundle method}
\label{sec: bundle_method}
We use an adaptation of a proximal bundle method utilised by \citet{kim2022scalable} to solve the \textit{p}-Lagrangian dual problem \eqref{eq: dual_problem}. The bundle method relies on an iterative approximation of the \textit{p}-Lagrangian dual function $L(\mu)$ with piece-wise linear functions via cutting planes. \citet{kim2022scalable} proposed a solution method for stochastic mixed-integer programming problems inspired by \citet{caroe1999dual} and utilised the adaptation of the proximal bundle method presented in \citep{kiwiel1990proximity} as a subroutine to solve dual problems. \citet{kiwiel1990proximity} developed a new approach for updating the weights of proximal terms in bundle methods for minimizing a convex function \citep{kc1984algorithm, kiwiel1985exact, kiwiel1987constraint}. This technique can significantly reduce the number of cutting planes required to reach the desired convergence accuracy. The pseudo-code of the adaption of the proximal bundle method proposed by \citet{kim2022scalable} is presented in Algorithm \ref{al: proximal bundle method}. 

Suppose that in the $k^\text{th}$ iteration of the proximal bundle method, we have computed Lagrangian multipliers $\mu_l$ and centres of mass $\overline{\mu}_l,$ for $l = 1, \dots, k-1$. In what follows, we present the adaptation of the proximal bundle method to update these parameters. 

The Lagrangian multiplier $\mu_k$ is computed as follows 
\begin{align}
   \mu_k = \argmax_\mu \left\{m_k(\mu) - \frac{u_k}{2} || \mu_k - \overline{\mu}_{k-1}||\right\}, \label{bm: cutting_plane_solution}
\end{align}
where $m_k(\mu)$ is piece-wise linear approximation of $L(\mu)$ at iteration $k$ given by 
\begin{align}
    m_k(\mu) = & \max_{\theta^s} \sum_{s \in S } \theta^s  \label{bm: cutting_plane_objective} \\ 
    & \text{s.t.: }  \theta^s \le L^s(\mu_l) - \left(\frac{\partial L^s(\mu_l)}{\partial\mu_l} \right)^\top (\mu - \mu_l), \  \forall s \in S, \ l = 1, \dots, k-1. \label{bm: cutting_plane_constraint}
\end{align}


The convergence of the proximal bundle method strongly relies on the update of the proximal parameter $u_k$ and of the centre of mass of $\overline{\mu}_k$. In line with the procedure the developed in \cite{kim2022scalable}, the centre of mass $\overline{\mu}_k$ is updated as follows
\begin{align}
    \overline{\mu}_{k} = \begin{cases} \mu_{k}, \ \ \ \ \text{if } L(\mu_k) \ge L(\overline{\mu}_k) + a_L v_k\ & \text{(serious step)}\\   \overline{\mu}_{k-1}, \  \text{otherwise} \ & \text{(null step)}, \end{cases} \label{bm: centre_update}
\end{align}
where we typically have $a_L \in (0, 0.5)$ and 
\begin{align}
    v_k = m_k(\mu_k) - L(\overline{\mu}_{k-1}) \label{bm: v_k}
\end{align}
representing the predicted increase of \textit{p}-Lagrangian function $L(\mu)$.

The proximal term $u_k$ must be chosen carefully. To prevent the proximal bundle method from taking a serious step too frequently (after too little improvement in $L(\mu)$), $u_k$ cannot be too large. On the other hand, if $u_k$ value is too small, the method will take many null steps before it finds a good candidate for the new centre of mass. To accelerate the performance of the proximal bundle method, tests identifying whether the proximal parameter $u_k$ value was too small or too large can be employed. 

The case when $u_k$ is too large can be identified by testing whether 
\begin{align} 
    L(\mu_k) \ge L(\overline{\mu}_{k-1}) + a_R v_k, \label{bm: too_large_test}
\end{align}
where $a_R \in (a_L, 1)$. If \eqref{bm: too_large_test} holds the proximal term $u_k$ is updated as 
\begin{align}
   u_{k+1} = \max\{ h_k, C^u_{min} u_k, u_{min}\}, \label{bm: too_large_update} 
\end{align}
with
\begin{align}
    h_k = 2 u_k \left( 1 - \frac{L(\mu_k) - L(\overline{\mu}_{k-1})}{v_k}\right), \label{bm: h_k}
\end{align}
and $C^u_{min} \in \mathbb{R}$. On the other hand, whether the proximal term $u_k$ is too small is identified by the test 
\begin{align}
    \overline{\delta}_k > \max \{ \delta_k(\overline{\mu}_{k-1}) + \lvert g_k \rvert, C^v v_k \}, \label{bm: too_small_test}
\end{align}
where ${C^v \in \mathbb{R}}$ and 
\begin{align}
    \overline{\delta}_k = L(\mu_k) -  \left(\sum_{s \in S}\frac{\partial L^s(\mu_k)}{\partial\mu_k} (x^s_k) \right)^\top (\mu_k - \overline{\mu}_{k-1}) - L(\overline{\mu}_{k-1}), \label{bm: overline_delta_k} 
\end{align}
in which $x^s_k$, $\forall s \in S$, is the optimal solution of the \textit{p}-Lagrangian sub-problem $L^s(\mu)$ with $\mu = \mu_k$,
\begin{align}
    g_k \in \partial m_k(\mu_k), \label{bm: g_k}
\end{align}
and 
\begin{align}
  \delta_k(\mu) = m_k(\mu_k) + (g_k)^\top(\mu -\mu_k) -L_{\pi}(\mu), \label{bm: delta_k}
\end{align}
where $\partial$ denotes the subdifferential of $m_k$ at $u_k$, making thus $g_k$ a subgradient of $m_k$ at $u_k$. If \eqref{bm: too_small_test} holds, the proximal term $u_k$ is updated as 
\begin{align}
   u_{k+1} = \min\{ h_k, C^u_{max} u_k\}, \label{bm: too_small_update}
\end{align}
where $C^u_{max} \in \mathbb{R}$. Algorithm \ref{al: proximal bundle method} summarises the developed proximal bundle method, starting with a step $k = 0$ and the initialisation of the parameters.

 \begin{algorithm}[h]
		\caption{Proximal bundle method}
		\label{al: proximal bundle method}
		\begin{algorithmic}
	    \State \textbf{initialise}: $k = 0, k_{max}, \epsilon_{BM}, \mu_0, \overline{\mu}_0 = \mu_0, u_{min}, u_1 > u_{min}$, $C^u_{min}$, $C^u_{avg}$, $C^u_{max}$, $C^v$, $i_{min}$, $i_{max}$ and $i^u_1 = 0$.
	    \State For each $s \in S$ solve $L^s(\mu)$ with $\mu = \mu_0$ and solve \eqref{bm: cutting_plane_objective}--\eqref{bm: cutting_plane_constraint} with $l = 0$ to form $m_1$.
	    \Repeat 
	    \State $k = k+1$.
	    \State From \eqref{bm: cutting_plane_solution} obtain $\mu_k$ and the value of $m_k(\mu_k)$.
	    \State For each $s \in S$, solve $L^s(\mu)$ at the point $m = m_k$.
	    \State Compute  $v_k, h_k, \overline{\delta}_k, g_k$ and $\delta_k$ as in \eqref{bm: v_k}, \eqref{bm: h_k}, \eqref{bm: overline_delta_k}, \eqref{bm: g_k} and \eqref{bm: delta_k}, respectively. 
	    \If{$ L(\mu_k) - L(\overline{\mu}_{k-1}) \ge a_L v_k $}
	    \State $\overline{\mu}_k =\mu_k$
	    \If{$L(\mu_k) - L(\overline{\mu}_{k-1}) \ge a_R v_k$ and $i^u_k > 0 $}
	    \State $u_{k+1} = \max\{h_k, C^u_{min}u_k, u_{min}\}$
	    \ElsIf{$i^u_k > i_{max}$}
	    \State $u_{k+1} = \max\{C^u_{avg}u_k, u_{min}\}$
	    \EndIf
	    \State $i^u_{k+1} = \begin{cases} \max\{i^u_k+1, 1\}, & \text{if } u_{k+1} = u_k  \\
	                           1, & \text{otherwise} \end{cases} $
	    \Else
	    \State $\overline{\mu}_k =  \overline{\mu}_{k-1}$
	    \If{$\overline{\delta}_k > \max \{ \delta_k + \lvert g_k \rvert, C^vv_k \}$ and $i^u_k < i_{min}$}
	    \State $u_{k+1} = \min\{h_k, C^u_{max}u_k\}$
	    \Else
	    \State $u_{k+1} = C^u_{max}u_k$
	    \EndIf
	    \State $i^u_{k+1} = \begin{cases} \min\{i^u_k-1, -1\}, & \text{if } u_{k+1} = u_k,  \\
	                           -1, & \text{otherwise}. \end{cases} $
	    \EndIf
	    \State Formulate $m_{k+1}$ as in \eqref{bm: cutting_plane_objective}--\eqref{bm: cutting_plane_constraint}.
	    \Until $v_k \le \epsilon_{BM}$ \textbf{or} $k > k_{max}$ 
	        
	    \State \textbf{return:} $\overline{\mu}_k$, $L(\overline{\mu}_k), \ (x^s_{k_{max}}, y^s_{k_{max}}, w^s_{k_{max}})_{s \in S}$ solving $L(\overline{\mu}_k)$.
		\end{algorithmic}%
\end{algorithm}%
 Following the developments in \cite{kim2022scalable}, algorithm \ref{al: proximal bundle method} includes an additional parameter $i^u_k$ that counts consecutive serious or null steps and enforces the tuning of the proximal term $u_k$, hoping to speed up the algorithm's convergence. The algorithm terminates when predicted increases $v_k$ are within an arbitrary tolerance $\epsilon_{BM}$. For proof of the convergence of the bundle method adaptation presented in Algorithm \ref{al: proximal bundle method}, one can refer to, for instance, \cite{kiwiel1990proximity}.

\subsubsection{Frank-Wolfe Progressive-Hedging method} 
\label{sec: FWPH}
Alternatively, one can apply the Frank-Wolfe Progressive-Hedging (FWPH) method \citep{FW_PH_2018} to solve the $p$-Lagrangian dual problem \eqref{eq: dual_problem}. FWPH is applied to the primal characterisation of \eqref{eq: dual_problem}: 
\begin{align}
    z_{LD} = \mini_{x, \overline{x}, y, w} \left\{ \begin{aligned} \sum_{s \in \mathcal{S}}\  \pi^s \left(c^\top x^s + q^{s^\top} y^s + \sum _{(i,j) \in B_Q} Q^s_{i,j} w^s_{i,j} \right) \\ : (x^s,y^s,\Gamma^s) \in \text{conv}(G^s), x^s = \overline{x}, \ \forall s \in S \end{aligned} \right\},  \label{primal_characterisation_dual_problem} 
\end{align}
where conv$(G^s)$ denotes the convex hull of $G^s$ for each $s \in S$. 

The FWPH method primarily relies on the classical progressive hedging method \citep{rockafellar1991scenarios}. However, unlike progressive hedging, FWPH can guarantee convergence in case the original problem is mixed-integer. Indeed, using the progressive hedging method as proposed in \cite{rockafellar1991scenarios} to solve 2SSMIP might result in suboptimal bounds, cycling behaviour and poor convergence behaviour of Lagrangian dual bound for problem \eqref{eq: dual_problem} as the presence of integer variables hinders its convergence guarantees. As a result, progressive hedging has typically been employed as a heuristics approach (see, for example, \cite{watson2011progressive}). FWPH integrates an extension of the Frank-Wolfe method called the simplicial decomposition method (SDM) to iteratively construct an inner approximation of conv$(G^s)$ for each $s \in S $. The composition of SDM and progressive hedging method allows for overcoming the aforementioned convergence issue. Additionally, it allows replacing the additional step of solving mixed-integer linear sub-problems with solving convex continuous quadratic sub-problems when calculating the Lagrangian dual bound. This, in turn, improves the computational performance of the FWPH method \citep{FW_PH_2018}.

The FWPH method uses the augmented Lagrangian dual problem, i.e., a modified Lagrangian dual problem in which the Lagrangian dual function is augmented by a penalty term that acts as a regularisation term. The augmented Lagrangian dual function based on relaxing the NAC constraints $x^s = \overline{x}, \forall s \in S$ in RNMDT$_p$ problem \eqref{RNMDT_p} is 
\begin{align}
    L_{\rho}(x, y, w, \overline{x}, \mu) = \sum_{s \in \mathcal{S}}\  \pi^s L^s_{\rho}(x^s, y^s, w^s, \overline{x}, \mu^s), \label{aug_p_LR}
\end{align}
where 
\begin{equation*}
 L^s_\rho(x^s, y^s, w^s, \overline{x}, \mu^s) = c^\top x^s + q^{s^\top} y^s + \sum _{(i,j) \in B_Q} Q^s_{i,j} w^s_{i,j} + {\mu^s}^\top (x^s - \overline{x}) + \frac{\rho}{2} {\|x^s - \overline{x} \|}^2_2
\end{equation*} 
and $\rho > 0$ is a penalty parameter. 

The FWPH algorithm pseudo-code is stated in Algorithm \ref{al: FWPH}. The parameter $k_{max}$
 defines the maximum number of iterations for the Frank-Wolfe method and $\epsilon_{FWPH}$ is a convergence tolerance parameter. The termination criterion involves the term $\sum_{s \in S} \pi^s {\left\lVert x^s_k - \overline{x}_{k-1} \right\rVert}$ that represents the sum of squared norms of primal and dual residuals associated with \eqref{primal_characterisation_dual_problem}. These residuals evaluate how close the solution candidate $((x^s,y^s,w^s), \overline{x})$ is to satisfy the necessary and sufficient optimality conditions for \eqref{primal_characterisation_dual_problem}.

 \begin{algorithm}[h]
		\caption{Frank-Wolfe progressive hedging (FWPH) method}
		\label{al: FWPH}
		\begin{algorithmic}
		\State \textbf{initialise}: $(V^s_0)_{s \in S},  (x^s_0)_{s \in S},  \mu_0, \rho, \alpha, \epsilon_{FWPH}, k_{max}$, $t_{max}$ and $\epsilon_{SDM}$.
		\State Compute  $\overline{x}_0 =\sum_{s \in S } \pi^s x^s_0 $ and  $\mu^s_1 =\mu ^s_0 + \rho(x^s_0 - \overline{x}_0)$.
		\For{ $k = 1, \dots, k_{max}$}
		
		   \For{$s \in S$}
		    
		        \State $\tilde{x}^s = (1-\alpha) \overline{x}_{k-1} + \alpha x^s_{k-1}$,
		        \State $[x^s_k, y^s_k, w^s_k, V^s_k, L^s(\mu_k^s)] =
		        SDM(V^s_{k-1}, \tilde{x}^s, \mu^s_k, \overline{x}_{k-1}, t_{max}, \epsilon_{SDM})$
		   \EndFor%
		   \State Compute $L(\mu_k) =\sum_{s \in S} \pi^s L^s(\mu^s_k)$ and $\overline{x}_k =\sum_{s \in S} \pi^s x^s_k$.
		   \If{$\sqrt{\sum_{s \in S} \pi^s {\left\lVert x^s_k - \overline{x}_{k-1} \right\rVert}^2_2}  \le \epsilon_{FWPH}  $}
		   \State \textbf{return} $((x^s_k, y^s_k, w^s_k)_{s \in S}, \ \overline{x}_k, \ \mu_k, \ L(\mu_k))$
		   \EndIf%
		   
		   \State Compute $\mu^s_{k+1} =\mu^s_k + \rho(x^s_k - \overline{x}_k)$  for each $s \in S$.
 		   
		\EndFor%
		
		\State \textbf{return} $ (x^s_{k_{max}}, y^s_{k_{max}}, w^s_{k_{max}})_{s \in S}, \overline{x}_{k_{max}}, \mu_{k_{max}}, L(\mu_{k_{max}})$.
		\end{algorithmic}%
\end{algorithm}%

 As a subroutine, Algorithm \ref{al: FWPH} employs the simplicial decomposition method (SDM)  to minimise \linebreak $L_\rho^s(x, y, w, \overline{x}, \mu^s)$ over $(x,y,w) \subset \text{conv}(G^s)$ for a given $s \in S$. The pseudo-code for SDM is stated in Algorithm \ref{SDM}. The precondition for the SDM algorithm is that $V^s_0 \subset$ conv$(G^s)$ and $\overline{x} = \sum_{s \in S} \pi^s x^s_0$, where $V^s_t$ are discrete sets of points such that  $V^s_t \subset$ conv$(G^s)$. Parameter $t_{max}$ defines the maximum number of iterations for SDM, and $\epsilon_{SDM} > 0 $ is the convergence tolerance. The parameter $\alpha $ affects the initial linearisation point $\tilde{x}^s$ of the SDM method.
 
 \begin{algorithm}[h]
		\caption{Simplicial decomposition method (SDM)}
		\label{SDM}
		\begin{algorithmic}
	    \State \textbf{initialise}: $V^s_0, x^s_0, \mu^s, \overline{x}, t_{max}$ and $\epsilon_{SDM}$.
	    \For{$t = 1, \dots, t_{max}$}
        \State $\hat{\mu^s_t} =\mu^s + \rho(x^s_{t-1} - \overline{x})$,
	    \State $(\hat{x}^s, \hat{y}^s, \hat{w}^s) \in {\argmin}_{x,y,w}\Big\{ (c + \hat{\mu}^s_t)^\top x + q^{s\top} y + \sum _{(i,j) \in B_Q} Q^s_{i,j} w_{i,j}:$ 
        \State $\ \ \ \ \ \ \ \ \ \ \ \ \ \ \ \ \ \ \ \ \ \ \ \ \ \ \ \ \ \ \ \ \ \ \ \ \ \ \ \ \ \ \ \ \ \ \ \ \ \ \ \ \ \ \ \ \ \ \ \ \ \ \ \ \ \ \ \ \ (x,y,w) \in G^s \Big\}$
	    \If{$t = 1$}
	       \State $L^s(\hat{\mu}^s_t) =(c + \hat{\mu}^s_t)^\top \hat{x}^s + q^{s\top} \hat{y}^s + \sum _{(i,j) \in B_Q} Q^s_{i,j} \hat{w}^s_{i,j}$
	        
	    \EndIf%
     \State Compute
	    \State $\Gamma^t =- \Big[(c + \hat{\mu}^s_t)^\top (\hat{x}^s - x^s_{t-1}) + q^{s\top} (\hat{y}^s - y^s_{t-1})$.
     \State $ \ \ \ \ \ \ \ \ \ \ \ \ \ \ \ \ \ \ \ \ \ \ \ \ \ \ \ \ \ \ \ \ \ \ \ \ \ \ \ \ \ \ \ \ \ \ \ \ \ + \sum _{(i,j) \in B_Q} Q^s_{i,j} (\hat{w}^s_{i,j} - w^s_{t-1, i,j})\Big]$,
	    \State $V^s_t =V^s_{t-1} \cup \{(\hat{x}^s, \hat{y}^s, \hat{w}^s)\}$ and 
     \State $(x^s_t, y^s_t) \in \argmin_{x,y,w}\left\{ L^s_\rho(x,y,w,\overline{x}, \hat{\mu}^s_t) : (x,y,w) \in \text{conv}(V^s_t)\right\}$.
	    \If{$\Gamma^t \le \epsilon_{SDM}$}
	        \State \textbf{return} $(x^s_t, y^s_t, w^s_t, V^s_t, L(\hat{\mu}_t))$
	    \EndIf%

	    \EndFor%
	    
	    \State \textbf{return} $(x^s_{t_{max}}, y^s_{t_{max}}, w^s_{t_{max}}, V^s_{t_{max}}, L(\hat{\mu}_{t_{max}}))$.
	    
		\end{algorithmic}%
\end{algorithm}%

\section{Dual decomposition}
\label{sec:dual_bnb}
In this section, we present the branching approach we employ, which is inspired by dual decomposition proposed in \cite{caroe1999dual}. The authors proposed a solution method for linear stochastic multi-stage problems that may involve integrality requirements at each stage. The solution method relies on dual decomposition combined with branch-and-bound strategies to ensure convergence. In what follows, we discuss our adaptation of the solution method proposed in \cite{caroe1999dual} for the mixed-integer RNMDT relaxations of RDEM problems.

Let $T$ be the set of unexplored nodes in the branch-and-bound search, in which each node is denoted by $N$. The key idea behind our approach is to extend the branch-and-bound procedure proposed in \cite{caroe1999dual} for the RNMDT$_p$ problem \eqref{RNMDT_p}. Specifically, we perform branching on the first-stage variables and use the solution of \textit{p}-Lagrangian dual problem, as described in \eqref{eq: dual_problem}, as the bounding procedure. To form candidates for feasible first-stage variables solution, the method uses an average $\overline{x}_N = \sum_{s \in S} \pi^s x_N^{*,s}$, combined with a rounding heuristic to fulfil the integrality requirements, where $x_N^{*,s}, \forall s \in S$, is obtained from solving the $N$ node-corresponding dual problem \eqref{eq: dual_problem}.

If $\overline{x}_N$ violates integrality conditions for some integer index $i$, i.e., $\lfloor \overline{x}_{N,i} \rfloor < \overline{x}_{N, i} < \lceil \overline{x}_{N,i} \rceil$, two nodes $N^L$ and $N^R$ with the correspondent sub-problems \eqref{eq: dual_problem} are formed from parent node $N$, where feasibility sets $G^s_{N^L}$ and $G^s_{N^R}$, $ \forall s \in S$, are formed respectively as
\begin{align}
 & G^s_{N^L} = G^s_{N} \cap \left\{ x^s_{i} \le \lfloor \overline{x}_{N,i} \rfloor \right\} \text{ and }
 \label{eq: int_left_branching} \\
 & G^s_{N^R} = G^s_{N} \cap \left\{ x^s_{i} \ge \lceil \overline{x}_{N,i} \rceil \right\}.
 \label{eq: int_right_branching}
\end{align}

There are multiple approaches for selecting the integer index $i$ of fractional-valued variable, i.e., $\lfloor \overline{x}_{N,i} \rfloor < \overline{x}_{N, i} < \lceil \overline{x}_{N,i} \rceil$ to perform branching. \cite{MORRISON201679} describes several common methods, of which we highlight the following:

\begin{itemize} \item Branching on the most (least) fractional or most (least) infeasible variable (see, for example, \cite{ACHTERBERG200542, ortega2003branch}): This approach involves selecting the integer index $i$ such that the fractional part of $\overline{x}_{N,i}$ is closest to (or furthest from) 0.5.

\item Pseudocost branching (see, for example, \cite{benichou1971experiments}): In this method, the integer index $i$ is chosen based on the expected impact of branching on $\overline{x}_{N,i}$ on the objective function, leveraging the history of previous branching decisions in the tree. 

\item Strong branching (see, for example, \cite{ACHTERBERG200542, achterberg2007constraint}): This technique recommends selecting the integer index $i$ such that branching on $\overline{x}_{N,i}$ results in the greatest change to the objective function. This is typically evaluated by solving LP relaxations for the child nodes generated by branching on $\overline{x}_{N,i}$. 
\end{itemize}

In addition to these strategies, the authors in \cite{MORRISON201679} highlight several other methods, including combining strong and pseudocost branching (see, for example, \cite{linderoth1999computational}), backdoor branching, which uses an auxiliary integer program to identify a subset of indices that reduce the search tree size (see, for example, \cite{10.1007/978-3-642-20807-2_15}), and information-theoretic branching, where indices are chosen to minimize uncertainty in the sub-problems (see, for example, \cite{GILPIN2011147}).

In our approach, a predominant feature is that we typically have a small search tree and that each node typically requires considerable computational work. The combination of these two characteristics led us to lean towards a simpler branching strategy.  

If $\overline{x}_N$ satisfies integrality conditions but $x_N^{*,s}, \ \forall s \in S$, violates non-anticipativity conditions, two nodes $N^L$ and $N^R$ with the correspondent sub-problems \eqref{eq: dual_problem} are formed from the parent node $N$, where feasibility sets $G^s_{N^L}$ and $G^s_{N^R}$, $\forall s \in S$, are formed respectively as
\begin{align}
 &G^s_{N^L} = G^s_{N} \cap \left\{ x^s_{i} \le \overline{x}_{N,i} - \epsilon_{BB} \right\} \text{ and }
  \label{eq: NAC_left_branching} \\
 &G^s_{N^R} = G^s_{N} \cap \left\{ x^s_{i} \ge \overline{x}_{N,i} +  \epsilon_{BB} \right\},
 \label{eq: NAC_right_branching}
\end{align} %
where $\epsilon_{BB} > 0$. The branching index $i$ is chosen based on the measure of the dispersion in the first-stage scenario solutions, e.g., if the dispersion of the component $i: \ \sigma_i = \max_{s \in S} x_{N,i}^{*,s} - \min_{s \in S} x_{N,i}^{*,s}$  is zero, this should imply the non-anticipativity of this component
$$
x^{*,1}_{N,i} = \dots = x^{*,\mid S \mid}_{N,i}.
$$ 
Therefore, in case of violating non-anticipativity constraints, branching is performed on the index $i$ with the largest dispersion. 

The branch-and-bound search requires that any node not pruned or fathomed must be explored. Therefore, reducing the size of the search tree is directly tied to the efficiency of the pruning strategy. \cite{MORRISON201679} highlight the following major groups of strategies
\begin{itemize}
\item Pruning by bound: This strategy involves pruning nodes whose sub-problem objective value is worse (higher in the case of minimisation) than the incumbent solution objective value (see, for example, \cite{arora2002proving, VILA2014105}). The bounds of the sub-problems can be calculated either by solving a linear relaxation of the node sub-problems or by using duality techniques.
\item Pruning by the dominance relations: this strategy involves pruning a node if another node’s sub-problem dominates its sub-problem. A sub-problem is said to be dominated if another sub-problem guarantees a better solution or an equally good solution that is less computationally intensive (see, for example, \cite{sewell_bbr_2012, fischetti2010pruning}). Dominance rules can be defined based on whether memory is used to store previously explored sub-problems. The \textit{memory-based dominance} involves comparing unexplored nodes to sub-problems of already generated and stored nodes. If a stored sub-problem dominates a new one, the new sub-problem can be pruned. The \textit{non-memory-based} dominance involves determining the existence of a dominating node without relying on stored sub-problems. It checks whether a dominance relation holds regardless of whether the dominating node has been explored. 
\end{itemize}
In the proposed approach, we prune nodes based on the lower bounds generated by solving the dual sub-problems, for the same reasons that we chose simpler branching rules (i.e., small search trees with computationally expensive nodes).

Algorithm \ref{al: dual_bnb} summarises adaptation of the branch-and-bound method presented by \cite{caroe1999dual} that hereafter we refer to as \textit{p}-BnB. For each branch-and-bound node $N \in T$, we generate node sub-problem \eqref{eq: dual_problem} and compute its dual bound value $z^*_N$ as well as corresponding optimal dual and primal variables values $(\mu_N^{*,s})_{s \in S}$ and  $(x_N^{*,s}, y_N^{*,s}, w_N^{*,s})_{s \in S}$, respectively, by applying Algorithm \ref{al: proximal bundle method} or \ref{al: FWPH}. If the dual bound value $z^*_N > z_{UB}$ or one of the $N$ sub-problems ${s \in S}$ is infeasible, the node $N$ is fathomed. Otherwise, we check whether solution $x_N^{*,s}$ violates non-anticipativity or integrality conditions. If that is the case, we perform branching as described in \eqref{eq: int_left_branching}--\eqref{eq: int_right_branching} on the most fractional variable $\overline{x}_{N,i}$ if $x_N^{*,s}$ violates integrality conditions. Otherwise, we perform branching as described in \eqref{eq: NAC_left_branching}--\eqref{eq: NAC_right_branching} on the variable with the largest dispersion $\sigma_i$ if $x_N^{*,s}$ violates non-anticipativity conditions. If $x_N^{*,s}$ satisfies both non-anticipativity and integrality conditions, we update the best upper bound value $z_{UB} = z^*_N$ and best solution value $x^* = \overline{x}_N$. Lastly, we update the best lower bound value $z_{LB}$ by setting it to the smallest dual bound value $z^*_N$ among the nodes $N$ that are yet to be fathomed. The algorithm continues until the set $T$ is empty. 

 \begin{algorithm}[h]
		\caption{\textit{p}-branch-and-bound method (\textit{p}-BnB)}
		\label{al: dual_bnb}
		\begin{algorithmic}
			\State \textbf{initialise}: $T =\emptyset $, $z_{UB} =\infty$, $z_{LB} =- \infty$, $x^* =\emptyset$, $\epsilon_{BB} > 0$ and $\epsilon_{NAC} \ge 0$.
			\State Create root node $N_0$ sub-problem \eqref{eq: dual_problem}, $T = T \cup \{N_0\} $.
		    \Repeat
		    \State Choose a node $N \in T$.
		    \State $T = T \setminus \{N\}$.
		    \State Apply Algorithm \ref{al: proximal bundle method} or \ref{al: FWPH} to the node $N$ sub-problem  \eqref{eq: dual_problem} to obtain $z_N^*$, $(\mu_N^{*,s})_{s \in S}$ and $(x_N^{*,s}, y_N^{*,s}, w_N^{*,s})_{s \in S}$.
                \If{$z^*_N > z_{UB}$ or one of the $N$ sub-problems is infeasible}  \State fathom $N$             
                \Else{}
    		    \State Compute $\overline{x}_N = \sum_{s \in \mathcal{S}} \pi^s x_N^{*,s}$.
    		    \State Compute $\sigma_i =\max_{s \in S} \braces{x_{N,i}^{*,s}} - \min_{s \in S} \braces{x_{N,i}^{*,s}}$ for $i \in \braces{1, \dots, n_x}$.
    		    \If{$\max_{i \in 1, \dots, n_x} \braces{\sigma_i} \le \epsilon_{NAC}$ }      
    		        \If{$\overline{x}_{N,i}$ is fractional for some integer index $ i \in \{1, \dots, n_x$\}} 
    		            \State Choose integer variable  index $ i \in  \braces{1, \dots, n_x}$ such as $\lfloor \overline{x}_{N,i} \rfloor < \overline{x}_{N,i} < \lceil \overline{x}_{N,i} \rceil$.  
    		            \State Create two new nodes $N^L$ and $N^R$ via \eqref{eq: int_left_branching} and \eqref{eq: int_right_branching}, respectively.
    		        \ElsIf{$z_{UB} > z_N^*$} 
    		            \State $z_{UB} = z_N^*$, 
    		            \State $x^* = \overline{x}_N$ 
    		        \EndIf
    		    \ElsIf{$\max_{i \in 1, \dots, n_x} \{\sigma_i\} > \epsilon_{NAC}$ and $z_{UB} > z_N^*$}
    		        \If{$\overline{x}_{N,i}$ is fractional for some integer index $ i \in     1, \dots, n_x$}
    		            \State Choose integer variable index $ i \in     1, \dots, n_x$ such as $\lfloor \overline{x}_{N,i} \rfloor < \overline{x}_{N,i} < \lceil \overline{x}_{N,i} \rceil$. 
    		            \State Create two new nodes $N^L$ and $N^R$ via \eqref{eq: int_left_branching} and \eqref{eq: int_right_branching}, respectively.
    		        \Else{}
    		            \State Choose continuous variable index $i \in \argmax_i \sigma_i$.
    		            \State Create two nodes $N^L$ and $N^R$ via \eqref{eq: NAC_left_branching} and \eqref{eq: NAC_right_branching}, respectively.
    		        \EndIf
    		        
    		        \State $T =T \cup \{N^L,N^R\} $.
    		   
    		    \EndIf  
                \EndIf
		    \State Update $Z_{LB}$.

		    \Until{$T = \emptyset $} 
			
		\end{algorithmic}%
\end{algorithm}%

\begin{figure}[h!]
    \centering
    \includegraphics[width=\textwidth]{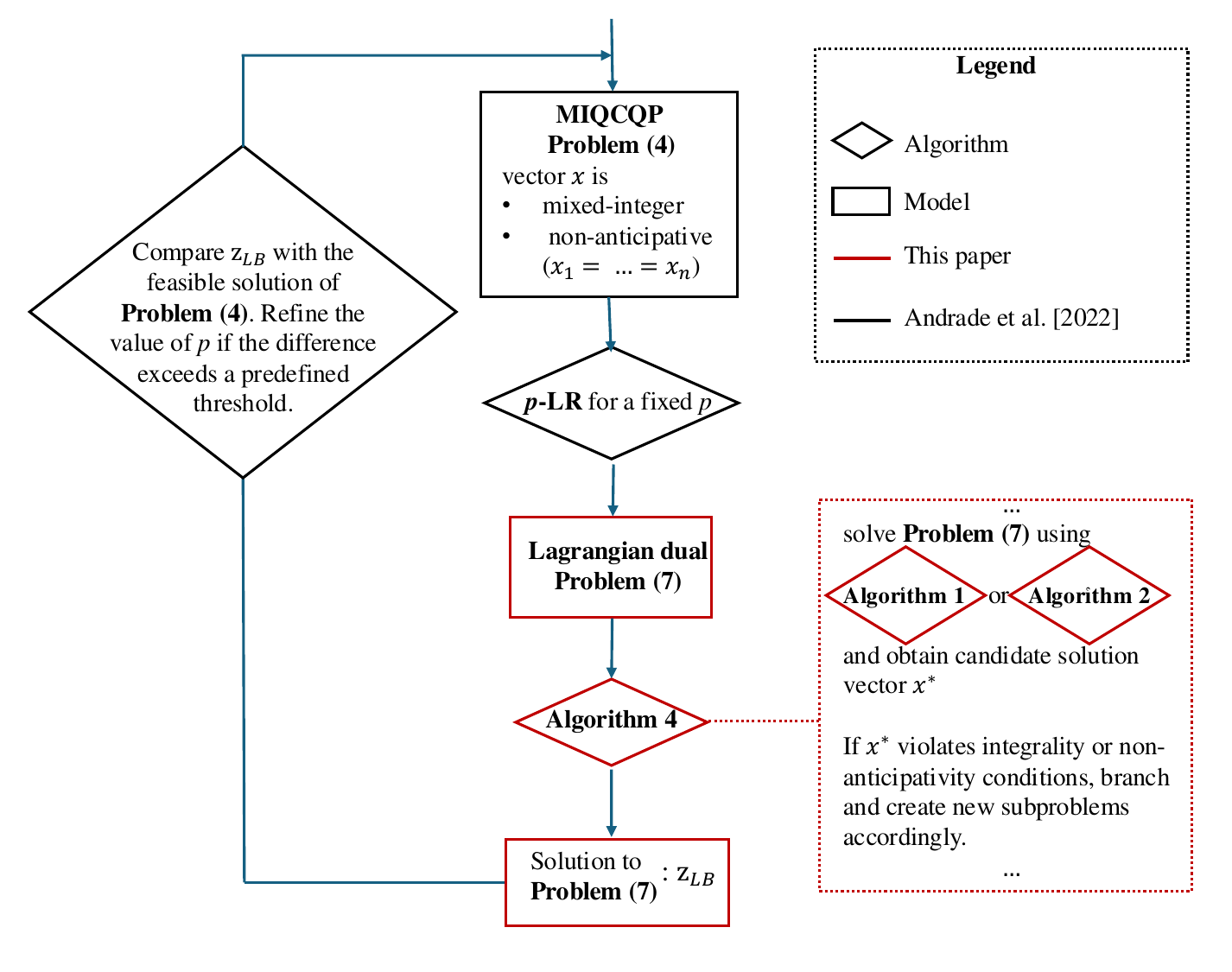}
    \caption{Graphical representation of Algorithm \ref{al: dual_bnb} structure}
    \label{fig: algorithm_structure}
\end{figure}

To clarify how the elements composing Algorithm \ref{al: dual_bnb} interact, we provide a graphical representation of the algorithm's structure in Figure \ref{fig: algorithm_structure}. Additionally, Section \ref{sec: ap3} presents an illustrative example, showcasing the procedures described in Algorithm \ref{al: dual_bnb} and demonstrating how they apply to a small problem instance.

In what follows, we provide a theoretical justification of the Algorithm \ref{al: dual_bnb} convergence to the optimal solution of RNMDT$_p$ relaxation (problem \eqref{RNMDT_p}). The convergence of the Algorithm \ref{al: dual_bnb} to the solution set of problem \eqref{RNMDT_p} considering any fixed value of $p = \{-\infty \dots, -1\}$ is stated in Theorem \ref{p_BNB_convergence}. Consequently, the solution set of problem \eqref{RNMDT_p} converges within a predefined tolerance level to the solution set of the primal RDEM (problem \eqref{eq: RDEM}) as the precision factor \textit{p} approaches $-\infty$. Formally, the justification for convergence of the RNMDT relaxation (problem \eqref{RNMDT_p}) is stated in Theorem \ref{RNMDT_convergence}. It is worth highlighting that when discussing the convergence from \eqref{RNMDT_p} to \eqref{eq: RDEM} as \textit{p} approaches $-\infty$ we take into consideration that, in practical applications, this corresponds to achieving a predetermined epsilon-accurate convergence rather than an absolute convergence.

\begin{theorem}
\label{p_BNB_convergence}
Suppose we consider the RNMDT relaxation (problem \eqref{RNMDT_p}) with an arbitrary fixed value of the precision factor $p = \{-\infty, \dots, -1\}$. Then Algorithm \ref{al: dual_bnb} converges to the solution $(x_N^{*,s}, y_N^{*,s}, w_N^{*,s})_{s \in S}$ that is optimal for problem \eqref{RNMDT_p}. 
\end{theorem}

\begin{prf}
In \cite{caroe1999dual}, the authors demonstrate the termination in finitely many steps and convergence of the Algorithm \ref{al: dual_bnb} to the optimal solution of problem \eqref{RNMDT_p} assuming that nodes \textit{p}-Lagrangian dual sub-problem \eqref{eq: dual_problem} are solved to optimality and hence, yielding optimal dual bound. Employing either Algorithm \ref{al: proximal bundle method} or \ref{al: FWPH} ensures the convergence to the optimal solution of the \textit{p}-Lagrangian dual sub-problem \eqref{eq: dual_problem}. For the convergence of Algorithms \ref{al: proximal bundle method} and \ref{al: FWPH} to the optimal solution of problem \eqref{p_LR}, please refer to the \cite{kim2022scalable} and \cite{FW_PH_2018}, respectively. 
\end{prf}


\begin{theorem}
\label{RNMDT_convergence}
Suppose we consider the RNMDT$_p$ relaxation problem \eqref{RNMDT_p} with an arbitrary fixed value of the precision factor $p = \{-\infty, \dots, -1 \}$. Then for any pair $(p_1, p_2)$ such that $p_1 < p_2 \le 0$  RNMDT$_{p_1}$ is a tighter (or equal) relaxation of the original RDEM problem than RNMDT$_{p_2}$.
\end{theorem}

\begin{prf}
See \cite[Theorem 6]{andrade2018rnmdt}.
\end{prf}




\section{Computational experiments}
\label{sec:numerical experiments}
This section presents numerical results for experiments performed using randomly generated non-convex 2SSMIP in the form of \eqref{eq: RDEM} or problems \eqref{eq: RDEM}, as we refer to them hereinafter. All code and instances generated are available on the GitHub repository \cite{belyak_p_bnb_2022}. The experiments were designed using Julia (Version 1.7.3) language \citep{Bazanson2017Julia} and Gurobi (Version 9.1) solver \citep{gurobi}. The code was run on Triton, Aalto University's high-performance computing cluster \citep{noauthor_triton_nodate}.

\subsection{Design of experiments}
We tested the efficiency of Algorithm \ref{al: dual_bnb} considering two alternative methods to solve sub-problems \eqref{eq: dual_problem}: the proximal bundle method (BM) presented in Section \ref{sec: bundle_method} and the Frank-Wolfe progressive hedging (FWPH) method presented in Section \ref{sec: FWPH}. Algorithm \ref{al: dual_bnb} was implemented using parallel computing, meaning that the scenario-sub-problems \eqref{eq: p_LR_decomposed_subproblem} are solved in parallel in both solution methods (proximal BM and FWPH). For each instance, the number of processes utilised for parallel computing was equal to 30. The computational efficiency of Algorithm \ref{al: dual_bnb} was compared with Gurobi's \cite{gurobi} branch-and-cut algorithm with standard parameterisation. 

We tested Algorithm \ref{al: dual_bnb} on 5 sets of randomly generated non-convex problem instances. Each set contained problems \eqref{eq: RDEM} with 50, 100 and 150 scenarios represented in two scales (small and large) as described in Table \ref{tab: parameters_low_density}. More detailed information on the range of the coefficients for problems \eqref{eq: RDEM} is presented in Appendix \ref{sec:ap1}. Additionally, we assumed two different values of the precision factor $p \in \{-2, -1\}$. Hence, we considered 60 instances in total. It is important to note that smaller values of $p$ were not explored, as the aforementioned values were found to be adequate for providing sufficient tolerance in terms of the RNMDT approximation. Furthermore, extensive numerical evidence demonstrating how good the RNMDT approximation is for different values of $p$ has already been extensively discussed in \cite{andrade2018rnmdt} and \cite{p_lagrangian_Andrade}.

For the sake of simplicity, for each instance, all the first-stage variables were assumed to be integer, and all the second-stage variables were assumed to be continuous. To make test instances similar to those available in library of test problems for stochastic integer programming \citep{stochastic_test_instances} (which are not MIQCQPs) in terms of the number of non-zero coefficients in the constraints and objective function, we assumed the quadratic matrices $Q^s$ and $U^s_m$ $\forall s \in S$, $\forall m \in M$ to be randomly generated with 1\% density, which is of the same order of the density (i.e., proportion of nonzero elements in the problems' constraint matrices) of the test problems considered.  

\begin{table}[h]

    \caption{Instance problems dimensions (per scenario)}
         \label{tab: parameters_low_density}
    \centering

    \begin{tabular}{ >{\centering}m{0.16\textwidth} >{\centering}m{0.25\textwidth}>{\centering}m{0.25\textwidth} c}
    \toprule
    \textbf{Instance size} & \textbf{ \# of $1^{\text{st}}$-stage variables} & \textbf{\# of $2^{\text{nd}}$-stage variables} & \textbf{\# of constraints} \\
    \midrule
    Small (S) & 100 & 100 & 100 \\ 
    Large (L) & 200 & 200 & 200 \\
    \bottomrule
    \end{tabular}
    
\end{table}

Therefore, the problems \eqref{eq: RDEM} with 50, 100 and 150 scenarios would have in total 5100, 10100 and 15100 variables, respectively, in the case of small (S) instances and 10200, 20200 and 30200 variables, respectively, in case of large (L) instances.

Table \ref{tab: algorithms_parametrisation} presents the parameter values used in the experiments for the proximal BM (Algorithm \ref{al: proximal bundle method}) and the FWPH (Algorithm \ref{al: FWPH}). In addition to the parameters stated in Table \ref{tab: algorithms_parametrisation}, the maximum number of iterations for the proximal BM  ($k_{max}$) was set to 1000. The maximum numbers of iterations for the FWPH algorithm ($k_{max}$) and simplicial decomposition method ($t_{max}$) were set to 1000 and 1, respectively. We considered $t_{max} = 1$ as \cite{boland2019parallelizable} suggest that considering higher values of $t_{max}$ parameter commonly increases the computational burden while not guaranteeing an adequate increase in the quality of the approximation of $\text{conv}(G^s)$ compared to the case when $t_{max} = 1$. It is also worth mentioning that since $t_{max} = 1$ for the simplicial decomposition method, the $\epsilon_{SDM}$ can be set to any arbitrary value as the condition $ \Gamma^t \le \epsilon_{SDM}$ is not going to be checked for the algorithm termination. Following \cite[Proposition 3.3]{FW_PH_2018}, since we assume $t_{max} = 1$, to guarantee the convergence of Algorithm \ref{al: FWPH} one should generate initial sets $\{(V_0^s)_{s \ in S}\}$ such that $\bigcap_{s \in S} \text{Proj}_{x}(\text{conv}(V_0^s)) \neq \emptyset$. Therefore, to initialise $(V^s_0)_{s \in S}$, we took one arbitrary scenario (in our case the first scenario in $S$, i.e., $s=1$) and set $V^1_0 = \{(x^1_0, y^1_0, w^1_0)\}$. Further, for each $s \in S, s \neq 1$, we initialised $V^s_0 = \{ {(x^s_0,y^s_0, w^s_0), (x^1_0, \overline{y}^s, \overline{w}^s)} \}$, where $(x^s_0, y^s_0, w^s_0)$ solves $L^s(\mu^s_0)$ and $(\overline{y}^s, \overline{w}^s))$ solves, 
$$
    \min_{y,w}\left\{  q^{s\top} y + \sum _{(i,j) \in B_Q} Q^s_{i,j} w_{i,j} : (x^1_0, y, w) \in G^s \right\}, \ \text{for each } \ s \in S.
$$ 

\begin{table}[h] 
     \caption{Algorithm parameters}
     \label{tab: algorithms_parametrisation}
    \centering
    \begin{tabular}{ l c }
        \toprule
        \multicolumn{2}{c}{\textbf{proximal bundle method}}\\
        \midrule
        $u_{min}$  &  $10^{-3}$ \\
        $m_{R}$  &  0.7 \\
        $m_{L}$  &  0.3 \\
        $i_{max}$ & 3\\
        $i_{min}$ & -3\\
        $C^u_{min}$ & 0.1\\
        $C^u_{avg}$ & 0.5\\
        $C^u_{max}$ & 10\\
        $C^v$ & 10\\
        $\epsilon_{BM}$ & $10^{-3}$\\ 
         \midrule
        \multicolumn{2}{c}{\textbf{Frank-Wolfe progressive hedging}}\\
        \midrule
         $\rho$  &  2 \\
         $\alpha$  &  1 \\
         $\epsilon_{FWPH}$ & $10^{-3}$ \\
         $\epsilon_{SDM}$ & not used \\
          \bottomrule
    \end{tabular} %
\end{table} %

 Starting dual multipliers values $\mu_0$ for Algorithms \ref{al: proximal bundle method} and \ref{al: FWPH} were set $\mu_0 = 0$. To set the first-stage variables $(x^s_0)_{s \in S}$ for Algorithm \ref{al: FWPH}, we considered the solution of the \textit{p}-Lagrangian dual function \eqref{p_LR} for a fixed value of the dual variable $\mu = \mu_0$. The tolerances $\epsilon_{BB}$ and $\epsilon_{NAC}$ for the \textit{p}-BnB (Algorithm \ref{al: dual_bnb}) were set to $10^{-6}$. As a time limit for solving each distinct instance, we considered one hour. In case multiple integer indices $i_1, \dots, i_{Int}$ are available for branching on the variable $\overline{x}_{N,i}$ in \textit{p}-BnB, we simply choose the first index $i_1$.

\subsection{Numerical results}
Table \ref{tab: numerical_results_small_density} presents averaged results of solving the small (S) and large (L) scale instances with the parameters as defined in Table \ref{tab: parameters_low_density} and $Q$ and $U$ matrices nonzero densities being set to 1\%. We compared the time required to solve the instances with the proposed \textit{p}-BnB method against solving them directly with the Gurobi solver (Full scale). The columns ``\textit{p}-BnB (FWPH)'' and ``\textit{p}-BnB (BM)'' report the solution for \textit{p}-BnB method when employing FWPH and proximal bundle method as a solution method for nodes sub-problems, respectively. Each cell in the ``Solution time'' section represents the average solution time value for 5 instances generated using 5 different random seeds. It is worth mentioning that when calculating the average value for the column ``Full-scale'' we have only considered the instances for which the Gurobi solver could generate a solution within one hour. 
\begin{table}[h]
\centering
\caption{Numerical results for the instances with low-density quadratic matrices}
\begin{tabular}{cccccc} \toprule
\multicolumn{3}{c}{\textbf{Instance parameters}}                                                                & \multicolumn{3}{c}{\textbf{Solution time (s)}}                                                                     \\ \hline
\multicolumn{1}{l}{Size}               & \multicolumn{1}{c}{$\absval{S}$} & \multicolumn{1}{c}{$p$} & \multicolumn{1}{l}{Full scale} & \multicolumn{1}{l}{\textit{p}-BnB (FWPH)} & \multicolumn{1}{l}{\textit{p}-BnB (BM)} \\ \hline
\multicolumn{1}{c}{\multirow{6}{*}{S}} & \multicolumn{1}{c}{50}    & -1                               & \multicolumn{1}{c}{83.88}      & \multicolumn{1}{c}{22.61}         & 15.68                          \\
\multicolumn{1}{c}{}                   & \multicolumn{1}{c}{50}    & -2                               & \multicolumn{1}{c}{131.22}     & \multicolumn{1}{c}{119.18}        & 10.18                          \\
\multicolumn{1}{c}{}                   & \multicolumn{1}{c}{100}   & -1                               & \multicolumn{1}{c}{185.56}     & \multicolumn{1}{c}{172.05}        & 41.01                          \\
\multicolumn{1}{c}{}                   & \multicolumn{1}{c}{100}   & -2                               & \multicolumn{1}{c}{358.34}     & \multicolumn{1}{c}{208.40}         & 55.36                          \\
\multicolumn{1}{c}{}                   & \multicolumn{1}{c}{150}   & -1                               & \multicolumn{1}{c}{316.23}     & \multicolumn{1}{c}{226.49}        & 50.28                          \\
\multicolumn{1}{c}{}                   & \multicolumn{1}{c}{150}   & -2                               & \multicolumn{1}{c}{535.56}     & \multicolumn{1}{c}{381.71}        & 92.61                          \\ \hline
\multicolumn{1}{c}{\multirow{6}{*}{L}} & \multicolumn{1}{c}{50}    & -1                               & \multicolumn{1}{c}{687.88}     & \multicolumn{1}{c}{630.49}        & 119.81                         \\
\multicolumn{1}{c}{}                   & \multicolumn{1}{c}{50}    & -2                               & \multicolumn{1}{c}{866.44}     & \multicolumn{1}{c}{420.63}        & 122.50                          \\
\multicolumn{1}{c}{}                   & \multicolumn{1}{c}{100}   & -1                               & \multicolumn{1}{c}{1505.92}    & \multicolumn{1}{c}{1637.15}       & 367.48                         \\
\multicolumn{1}{c}{}                   & \multicolumn{1}{c}{100}   & -2                               & \multicolumn{1}{c}{2490.45}    & \multicolumn{1}{c}{1708.13}       & 284.36                         \\
\multicolumn{1}{c}{}                   & \multicolumn{1}{c}{150}   & -1                               & \multicolumn{1}{c}{2463.96}    & \multicolumn{1}{c}{1372.53}       & 523.98                         \\
\multicolumn{1}{c}{}                   & \multicolumn{1}{c}{150}   & -2                               & \multicolumn{1}{c}{3412.82}    & \multicolumn{1}{c}{1031.00}          & 369.48   \\ \hline                     
\end{tabular}
\label{tab: numerical_results_small_density}
\end{table}

As the numerical results in Table \ref{tab: numerical_results_small_density} suggest, for small-scale instances, the proposed \textit{p}-BnB method outperformed commercial solver Gurobi in terms of the solution time regardless of the method employed to solve the dual sub-problems. This conclusion also applies to the large-scale instances, except for the instance with 100 scenarios and precision factor $p = -1$. On average, applying \textit{p}-BnB with Frank-Wolfe progressing hedging allowed for saving up 31.41\%  and 32.76\% of solution time for small- and large-scale instances, respectively, compared to solving the full-scale instances with Gurobi. The best improvement for the small-scale instances has been achieved for the instance with 50 scenarios and RNMDT precision factor $p=-1$, demonstrating a decrease in computational time by 73.05\% compared to solving the instance directly with Gurobi. For the large-scale instances, the largest reduction in solution time was observed for the instance with 150 scenarios and RNMDT precision factor $p = -2$, allowing for reducing the solution time required by Gurobi by 69.79\%. However, using \textit{p}-BnB with the proximal BM instead has demonstrated even further improvements in computational solution time. Compared to solving the full-scale instances with Gurobi, \textit{p}-BnB with the proximal BM demonstrated, on average, a decrease by 83.80\% and 83.42\% in solution time for the small- and large-scale instances, respectively. Moreover, the results suggest that the solution time improvement reached up to 92.24\% for the small-scale instances, as in the case of the instances with 50 scenarios and an RNMDT precision factor of $p = -2$. For the large-scale instances, the maximum improvement in solution time was achieved for the instance with 150 scenarios and RNMDT precision factor being $p = -2$, allowing a reduction of 89.17\% in the time required to solve that instance by Gurobi. However, It is important to highlight that in the case of utilising FWPH in the context of \textit{p}-BnB we have observed a considerable portion of computational time spent by FWPH on generating the sets $\{(V_0^s)_{s \ in S}\}$ at the beginning of Algorithm \ref{al: FWPH}, particularly for the instances with a high number of scenarios. Hence, addressing this issue could potentially lead to an improvement in its computational time performance, and possibly to approaching the efficiency of proximal BM.

Nevertheless, in all 60 instances, the \textit{p}-BnB explored only one (root) node to identify the optimal solution. This effectively means that all of these instances were such that there was no duality gap when solving the \textit{p}-Lagrangian duals and that bounds obtained by both methods were tight enough to find the optimal solution at the root node. This effect was also observed in \cite{FW_PH_2018} where the authors reported convergence of the FWPH method to the optimal solution for most of the stochastic mixed-integer problem instances. Additionally, the usage of \textit{p}-Lagrangian relaxation exploits the block-angular structure of the primal RDEM problem allowing one to obtain tighter bounds at the root if compared to linear-programming (LP) relaxation. Such phenomena have been reported in \cite{caroe1999dual} where the authors would obtain at a root node a duality gap of only 0.2– 0.3\% in case Lagrangian relaxation is explored while the LP-relaxation, however, would provide a duality gap of 2.0 –2.1\%

To demonstrate the convergence of the method in cases when the solution for the root node violates integrality or non-anticipativity conditions, we conducted another batch of experiments for somewhat less realistic instances in which the matrices $Q$ and $U$ densities are 90 \%. The increased densities of the matrices $Q$ and $U$ indicate a prevalence of non-zero coefficients of the bi-linear terms in the primal RDEM (problem \eqref{eq: RDEM}). This, in turn, implies that associated RNMDT relaxations (problem \eqref{RNMDT_p}) would have a significant number of auxiliary binary variables, thereby potentially increasing the likelihood for the existence of duality gap when solving corresponding \textit{p}-Lagrangian dual problems. However, to ensure convergence of \textit{p}-BnB within one hour, we tested \textit{p}-BnB on 5 instances with RNMDT precision factor $p = -1$ only and remaining parameters as before. Table \ref{tab: numerical_results_large_density} demonstrates the results of solving instances 1-5 with the proposed \textit{p}-BnB method employing the FWPH (\textit{p}-BnB (FWPH)) and proximal BM (\textit{p}-BnB (BM)) as a subroutine. The column “sol. time” reports the time required by Algorithm \ref{al: dual_bnb} to converge to an optimal solution with a 0.00\% gap, calculated as the relative difference between the upper bound (UB) and lower bound (LB) for the objective function generated by the corresponding method. The difference was calculated as $100 \% \frac{\text{UB} - \text{LB}}{\text{LB}} $. It is worth highlighting that solving full-scale instances with Gurobi resulted in convergence within one hour only for Instance 1, taking in a total of 2064.55 seconds. 

Interestingly, for the setting where the matrices are denser, FWPH outperforms BM as the Lagrangian dual solver in Instances 4 and 5, which is a consequence of FWPH leading the algorithm to explore less nodes than BM. In either case,  the number of nodes explored is small. As can be seen in Table \ref{tab: numerical_results_large_density}, the maximum number of nodes explored by \textit{p}-BnB was only 11, for Instance 5 using the proximal BM, while the average number of nodes explored was five. This is because despite the very high density of the quadratic matrices (90\%) in the instances, at the very first node, \textit{p}-BnB was able to generate a solution with a tight dual bound, on average being 0.01\%. In comparison, the average dual bound generated within one hour by solving the Instances 1-5 with Gurobi was 4.98\%. Further details on each individual Instance are provided in Appendix \ref{sec:ap2}.

\begin{table}[h]
    \caption{Dimensions of instances with high-density $Q$ matrices}
    \label{tab: parameters_high_density}
    \centering
    \resizebox{\textwidth}{!}{
    \begin{tabular}{ >{\centering}m{0.16\textwidth}>{\centering}m{0.16\textwidth} >{\centering}m{0.20\textwidth}>{\centering}m{0.20\textwidth} c}
    \toprule
    \textbf{Instance} & \textbf{\# of scenarios} & \textbf{\# of $1^{\text{st}}$-stage variables} & \textbf{\# of $2^{\text{nd}}$-stage variables} & \textbf{\# of constraints}   \\
    \midrule
    \textbf{1} & 15 & 30 & 25 & 25  \\ 
    \textbf{2} & 20 & 30 & 30 & 20 \\
    \textbf{3} & 20	& 40 & 15 & 15 \\
    \textbf{4} & 30	& 30 & 20 & 15 \\ 
    \textbf{5} & 40 & 20 & 10 & 15 \\
    \bottomrule
    \end{tabular}
    }
\end{table}

\begin{table}[h]
\centering
\caption{Numerical results for the instances with high-density $Q$ and $U$ matrices}
 \label{tab: numerical_results_large_density}

\begin{tabular}{ccccccc} \toprule
\multirow{2}{*}{\textbf{Instance}}  & \multicolumn{3}{c}{\textbf{\textit{p}-BnB (FWPH)}} & \multicolumn{3}{c}{\textbf{\textit{p}-BnB (BM)}} \\
                                  & sol. time (s)     & \# nodes    & \# iter.   & sol. time (s)   & \# nodes  & \# iter.  \\ \hline
{\bf 1}                              & 459.90   & 5           & 68      & 	114.42   & 1         & 30       \\
{\bf 2}                               &410.63   & 3           & 21         & 170.53  & 1         & 26        \\
{\bf 3}                               & 520.47	    & 5          & 145        & 925.73   & 3         & 312       \\
{\bf 4}                                & 374.79	    & 3           & 56        & 1439.22   & 5         & 268       \\
{\bf 5}                                & 427.94   & 9           & 191       &	3001.86   & 11         & 1525       \\ \bottomrule
\end{tabular}
\end{table}

\section{Conclusions} \label{sec:conclusion}

In this paper, we propose a novel method for solving two-stage stochastic programming problems whose deterministic equivalents are represented by non-convex MIQCQP models. Additionally, we assess the efficiency of this method by considering two alternative algorithms for solving dual sub-problems. The proposed method is named \textit{p}-branch-and-bound (\textit{p}-BnB) and combines a branch-and-bound-based algorithm inspired by \cite{caroe1999dual} with the \textit{p}-Lagrangian decomposition proposed in \cite{p_lagrangian_Andrade}. The \textit{p}-Lagrangian decomposition method relies on the composition of the mixed-integer-based relaxation of the non-convex MIQCQP problem using the reformulated normalized multiparametric disaggregation technique (RNMDT) \citep{andrade2018rnmdt} and classic Lagrangian relaxation. The construction of a mixed-integer-based relaxation for the primal non-convex MIQCQP problem is essential to ensure the validity of the dual bound associated with the primal problem. The \textit{p}-Lagrangian decomposition has been demonstrated to outperform the commercial solver Gurobi in terms of computational time required to generate the dual bounds for a primal non-convex MIQCQP problem, whose precision can be controlled by choice of parameters in RNMDT relaxation. However, \textit{p}-Lagrangian decomposition could not tackle the duality gap arising from the mixed-integer nature of the primal non-convex MIQCQP problems. In contrast, the proposed \textit{p}-BnB mitigates this issue by ensuring the integrality conditions of the optimal solution via a classic branch-and-bound approach. Additionally, following \cite{caroe1999dual}, the branch-and-bound procedure takes place whenever the first-stage variables candidates violate the non-anticipativity constraints. We also evaluated the efficiency of \textit{p}-BnB by considering two alternative solution methods for dual sub-problems in contrast to employing the classic bundle method as in the \textit{p}-Lagrangian decomposition \citep{p_lagrangian_Andrade}. We utilised Frank-Wolfe progressive hedging \citep{FW_PH_2018} and proximal bundle method \citep{kim2022scalable} to solve the node sub-problems. The Frank-Wolfe progressive hedging method is the enhancement of classic progressive hedging \citep{rockafellar1991scenarios} that guarantees convergence even for mixed-integer problems. The proximal bundle method is an enhancement of a classic bundle method \citep{lemarechal1974algorithm, zhao2002new} involving a new approach for updating the weights of proximal terms that can significantly reduce the number of iterations in the bundle method necessary to reach the desired convergence tolerance \citep{kiwiel1990proximity}. It is important to highlight that \textit{p}-BnB method ensures convergence to the solution of the non-convex MIQCQP problem given a predefined tolerance, as opposed to absolute convergence. Nonetheless, the former is generally sufficient for the majority of practical applications.

The \textit{p}-BnB efficiency has been tested on a set of RNMDT relaxations of randomly generated non-convex MIQCQP instances. Numerical experiments demonstrated the superior performance of the proposed \textit{p}-BnB method over attempts to solve full-scale RNMDT problems with the commercial solver Gurobi. Depending on the method utilised to solve dual sub-problems, the use of \textit{p}-BnB allowed for saving on average about 32 \% of the time required by Gurobi to solve RNMDT problem in case \textit{p}-BnB used Frank-Wolfe progressive hedging as a subroutine or about 84 \% of the time if proximal bundle method has been used. These results lead us to conclude that, although FWPH does potentially provide better dual bounds, it pays the price for being a more computationally demanding solution method. This trade-off flips when we consider artificially denser instances, in which, we can see that the FWPH method leads to fewer nodes explored in instances 4 and 5, and thus potential computational savings.

It is worth highlighting that the \textit{p}-BnB method implementation involves software engineering decisions that can greatly influence its performance. Nevertheless, our implementation still serves as a reliable proof of concept. Additionally, the \textit{p}-BnB method as proposed only considers rudimentary heuristics to generate feasible solutions for the primal RNMDT relaxation and the implementation of more sophisticated heuristics would likely improve the performance of \textit{p}-BnB, in a similar fashion as they are beneficial in mixed-integer programming solvers. Hence, one could further enhance \textit{p}-BnB computational efficiency. In particular, one potential path for improvement involves enhancing the branching strategies considered \citep{cornuejols2011improved, CHEN20031925}. Another possible direction could be an improvement of the FWPH method implementation. Additionally, an improvement of the procedure for generating the sets $\{(V_0^s)_{s \ in S}\}$ at the beginning of Algorithm \ref{al: FWPH} could bring new insight into \textit{p}-BnB performance and convergence rate when using FWPH as a subroutine.

\bibliographystyle{plainnat} 
\bibliography{references.bib}

\begin{thebibliography}{64}
\providecommand{\natexlab}[1]{#1}
\providecommand{\url}[1]{\texttt{#1}}
\expandafter\ifx\csname urlstyle\endcsname\relax
  \providecommand{\doi}[1]{doi: #1}\else
  \providecommand{\doi}{doi: \begingroup \urlstyle{rm}\Url}\fi

\bibitem[{Aalto scientific computing}(2022)]{noauthor_triton_nodate}
{Aalto scientific computing}.
\newblock Triton cluster, 2022.
\newblock URL \url{https://scicomp.aalto.fi/triton/#overview}.

\bibitem[Achterberg(2007)]{achterberg2007constraint}
Tobias Achterberg.
\newblock Constraint integer programming.
\newblock 2007.

\bibitem[Achterberg et~al.(2005)Achterberg, Koch, and Martin]{ACHTERBERG200542}
Tobias Achterberg, Thorsten Koch, and Alexander Martin.
\newblock Branching rules revisited.
\newblock \emph{Operations Research Letters}, 33\penalty0 (1):\penalty0 42--54, 2005.
\newblock ISSN 0167-6377.
\newblock \doi{https://doi.org/10.1016/j.orl.2004.04.002}.
\newblock URL \url{https://www.sciencedirect.com/science/article/pii/S0167637704000501}.

\bibitem[Ahmed et~al.(2015)Ahmed, Garcia, Kong, Ntaimo, Parija, Qiu, and Sen]{stochastic_test_instances}
S.~Ahmed, R.~Garcia, N.~Kong, L.~Ntaimo, G.~Parija, F.~Qiu, and S.~Sen.
\newblock Siplib: A stochastic integer programming test problem library, 2015.
\newblock URL \url{https://www2.isye.gatech.edu/~sahmed/siplib}.

\bibitem[Amos et~al.(1997)Amos, R{\"o}nnqvist, and Gill]{amos1997modelling}
F~Amos, M~R{\"o}nnqvist, and G~Gill.
\newblock Modelling the pooling problem at the new zealand refining company.
\newblock \emph{Journal of the Operational Research Society}, 48\penalty0 (8):\penalty0 767--778, 1997.

\bibitem[Andiappan(2017)]{andiappan_2017}
Viknesh Andiappan.
\newblock State-of-the-art review of mathematical optimisation approachesfor synthesis of energy systems.
\newblock \emph{Process Integration and Optimization for Sustainability}, 1\penalty0 (3):\penalty0 165--188, 2017.
\newblock \doi{10.1007/s41660-017-0013-2}.
\newblock URL \url{https://doi.org/10.1007/s41660-017-0013-2}.

\bibitem[Andrade et~al.(2019)Andrade, Oliveira, Hamacher, and Eberhard]{andrade2018rnmdt}
Tiago Andrade, Fabricio Oliveira, Silvio Hamacher, and Andrew Eberhard.
\newblock Enhancing the normalized multiparametric disaggregation technique for mixed-integer quadratic programming.
\newblock \emph{Journal of Global Optimization}, 73\penalty0 (4):\penalty0 701--722, April 2019.
\newblock ISSN 1573-2916.
\newblock \doi{10.1007/s10898-018-0728-9}.
\newblock URL \url{https://doi.org/10.1007/s10898-018-0728-9}.

\bibitem[Andrade et~al.(2022)Andrade, Belyak, Eberhard, Hamacher, and Oliveira]{p_lagrangian_Andrade}
Tiago Andrade, Nikita Belyak, Andrew Eberhard, Silvio Hamacher, and Fabricio Oliveira.
\newblock The p-lagrangian relaxation for separable nonconvex miqcqp problems.
\newblock \emph{Journal of Global Optimization}, 84\penalty0 (1):\penalty0 43--76, 2022.

\bibitem[Arora et~al.(2002)Arora, Bollob{\'a}s, and Lov{\'a}sz]{arora2002proving}
Sanjeev Arora, B{\'e}la Bollob{\'a}s, and L{\'a}szl{\'o} Lov{\'a}sz.
\newblock Proving integrality gaps without knowing the linear program.
\newblock In \emph{The 43rd Annual IEEE Symposium on Foundations of Computer Science, 2002. Proceedings.}, pages 313--322. IEEE, 2002.

\bibitem[Bashiri et~al.(2021)Bashiri, Nikzad, Eberhard, Hearne, and Oliveira]{bashiri2021two}
Mahdi Bashiri, Erfaneh Nikzad, Andrew Eberhard, John Hearne, and Fabricio Oliveira.
\newblock A two stage stochastic programming for asset protection routing and a solution algorithm based on the progressive hedging algorithm.
\newblock \emph{Omega}, 104:\penalty0 102480, 2021.

\bibitem[Belotti(2023)]{couenne}
Pietro Belotti.
\newblock Couenne: a user’s manual, 2023.
\newblock URL \url{https://www.coin-or.org}.

\bibitem[Belyak(2022)]{belyak_p_bnb_2022}
Nikita Belyak.
\newblock p-branch-and-bound metho.
\newblock \url{https://github.com/gamma-opt/p-BnB}, 2022.

\bibitem[B{\'e}nichou et~al.(1971)B{\'e}nichou, Gauthier, Girodet, Hentges, Ribi{\`e}re, and Vincent]{benichou1971experiments}
Michel B{\'e}nichou, Jean-Michel Gauthier, Paul Girodet, Gerard Hentges, Gerard Ribi{\`e}re, and Olivier Vincent.
\newblock Experiments in mixed-integer linear programming.
\newblock \emph{Mathematical programming}, 1:\penalty0 76--94, 1971.

\bibitem[Bergamini et~al.(2005)Bergamini, Aguirre, and Grossmann]{BERGAMINI20051914}
Maria~Lorena Bergamini, Pio Aguirre, and Ignacio Grossmann.
\newblock Logic-based outer approximation for globally optimal synthesis of process networks.
\newblock \emph{Computers \& Chemical Engineering}, 29\penalty0 (9):\penalty0 1914--1933, 2005.
\newblock ISSN 0098-1354.
\newblock \doi{https://doi.org/10.1016/j.compchemeng.2005.04.003}.
\newblock URL \url{https://www.sciencedirect.com/science/article/pii/S0098135405001006}.

\bibitem[Berthold et~al.(2012)Berthold, Heinz, and Vigerske]{berthold2012extending}
Timo Berthold, Stefan Heinz, and Stefan Vigerske.
\newblock Extending a cip framework to solve miqcp s.
\newblock In \emph{Mixed integer nonlinear programming}, pages 427--444. Springer, 2012.

\bibitem[Bezanson et~al.(2017)Bezanson, Edelman, Karpinski, and Shah]{Bazanson2017Julia}
Jeff Bezanson, Alan Edelman, Stefan Karpinski, and Viral~B. Shah.
\newblock Julia: {A} {Fresh} {Approach} to {Numerical} {Computing}.
\newblock \emph{SIAM Review}, 59\penalty0 (1):\penalty0 65--98, January 2017.
\newblock ISSN 0036-1445, 1095-7200.
\newblock \doi{10.1137/141000671}.
\newblock URL \url{https://epubs.siam.org/doi/10.1137/141000671}.

\bibitem[Boland et~al.(2018)Boland, Christiansen, Dandurand, Eberhard, Linderoth, Luedtke, and Oliveira]{FW_PH_2018}
Natashia Boland, Jeffrey Christiansen, Brian Dandurand, Andrew Eberhard, Jeff Linderoth, James Luedtke, and Fabricio Oliveira.
\newblock Combining progressive hedging with a frank--wolfe method to compute lagrangian dual bounds in stochastic mixed-integer programming.
\newblock \emph{SIAM JOURNAL ON OPTIMIZATION}, 28\penalty0 (2):\penalty0 1312--1336, 2018.
\newblock ISSN 1052-6234.
\newblock \doi{10.1137/16M1076290}.

\bibitem[Boland et~al.(2019)Boland, Christiansen, Dandurand, Eberhard, and Oliveira]{boland2019parallelizable}
Natashia Boland, Jeffrey Christiansen, Brian Dandurand, Andrew Eberhard, and Fabricio Oliveira.
\newblock A parallelizable augmented lagrangian method applied to large-scale non-convex-constrained optimization problems.
\newblock \emph{Mathematical Programming}, 175:\penalty0 503--536, 2019.
\newblock \doi{10.1007/s10107-018-1253-9}.

\bibitem[Bragalli et~al.(2012)Bragalli, D’Ambrosio, Lee, Lodi, and Toth]{bragalli2012optimal}
Cristiana Bragalli, Claudia D’Ambrosio, Jon Lee, Andrea Lodi, and Paolo Toth.
\newblock On the optimal design of water distribution networks: a practical minlp approach.
\newblock \emph{Optimization and Engineering}, 13\penalty0 (2):\penalty0 219--246, 2012.

\bibitem[Bush and Mosteller(1953)]{bush1953stochastic}
Robert~R Bush and Frederick Mosteller.
\newblock A stochastic model with applications to learning.
\newblock \emph{The Annals of Mathematical Statistics}, pages 559--585, 1953.

\bibitem[Calvo et~al.(2004)Calvo, de~Luigi, Haastrup, and Maniezzo]{calvo2004distributed}
Roberto~Wolfler Calvo, Fabio de~Luigi, Palle Haastrup, and Vittorio Maniezzo.
\newblock A distributed geographic information system for the daily car pooling problem.
\newblock \emph{Computers \& Operations Research}, 31\penalty0 (13):\penalty0 2263--2278, 2004.

\bibitem[Carøe and Schultz(1999)]{caroe1999dual}
Claus~C. Carøe and Rüdiger Schultz.
\newblock Dual decomposition in stochastic integer programming.
\newblock \emph{Operations Research Letters}, 24\penalty0 (1-2):\penalty0 37--45, February 1999.
\newblock ISSN 01676377.
\newblock \doi{10.1016/S0167-6377(98)00050-9}.
\newblock URL \url{https://linkinghub.elsevier.com/retrieve/pii/S0167637798000509}.

\bibitem[Castillo et~al.(2005)Castillo, Westerlund, Emet, and Westerlund]{castillo2005optimization}
Ignacio Castillo, Joakim Westerlund, Stefan Emet, and Tapio Westerlund.
\newblock Optimization of block layout design problems with unequal areas: A comparison of milp and minlp optimization methods.
\newblock \emph{Computers \& Chemical Engineering}, 30\penalty0 (1):\penalty0 54--69, 2005.

\bibitem[Castro(2016{\natexlab{a}})]{castro2016normalized}
Pedro~M Castro.
\newblock Normalized multiparametric disaggregation: an efficient relaxation for mixed-integer bilinear problems.
\newblock \emph{Journal of Global Optimization}, 64\penalty0 (4):\penalty0 765--784, 2016{\natexlab{a}}.

\bibitem[Castro(2016{\natexlab{b}})]{castro2016spatial}
Pedro~M Castro.
\newblock Spatial branch and bound algorithm for the global optimization of miqcps.
\newblock In \emph{Computer Aided Chemical Engineering}, volume~38, pages 523--528. Elsevier, 2016{\natexlab{b}}.

\bibitem[Cornu{\'e}jols et~al.(2011)Cornu{\'e}jols, Liberti, and Nannicini]{cornuejols2011improved}
Gerard Cornu{\'e}jols, Leo Liberti, and Giacomo Nannicini.
\newblock Improved strategies for branching on general disjunctions.
\newblock \emph{Mathematical Programming}, 130\penalty0 (2):\penalty0 225--247, 2011.

\bibitem[Cui et~al.(2013)Cui, Zheng, Zhu, and Sun]{cui2013convex}
XT~Cui, XJ~Zheng, SS~Zhu, and XL~Sun.
\newblock Convex relaxations and miqcqp reformulations for a class of cardinality-constrained portfolio selection problems.
\newblock \emph{Journal of Global Optimization}, 56\penalty0 (4):\penalty0 1409--1423, 2013.

\bibitem[Ding et~al.(2014)Ding, Bo, Li, and Sun]{ding2014bi}
Tao Ding, Rui Bo, Fangxing Li, and Hongbin Sun.
\newblock A bi-level branch and bound method for economic dispatch with disjoint prohibited zones considering network losses.
\newblock \emph{IEEE Transactions on Power Systems}, 30\penalty0 (6):\penalty0 2841--2855, 2014.

\bibitem[Feltenmark and Kiwiel(2000)]{feltenmark2000dual}
Stefan Feltenmark and Krzysztof~C Kiwiel.
\newblock Dual applications of proximal bundle methods, including lagrangian relaxation of nonconvex problems.
\newblock \emph{SIAM Journal on Optimization}, 10\penalty0 (3):\penalty0 697--721, 2000.

\bibitem[Fischetti and Monaci(2011)]{10.1007/978-3-642-20807-2_15}
Matteo Fischetti and Michele Monaci.
\newblock Backdoor branching.
\newblock In Oktay G{\"u}nl{\"u}k and Gerhard~J. Woeginger, editors, \emph{Integer Programming and Combinatoral Optimization}, pages 183--191, Berlin, Heidelberg, 2011. Springer Berlin Heidelberg.
\newblock ISBN 978-3-642-20807-2.

\bibitem[Fischetti and Salvagnin(2010)]{fischetti2010pruning}
Matteo Fischetti and Domenico Salvagnin.
\newblock Pruning moves.
\newblock \emph{INFORMS Journal on Computing}, 22\penalty0 (1):\penalty0 108--119, 2010.

\bibitem[Forrest et~al.(1974)Forrest, Hirst, and Tomlin]{forrest_practical_1974}
J.~J.~H. Forrest, J.~P.~H. Hirst, and J.~A. Tomlin.
\newblock Practical {Solution} of {Large} {Mixed} {Integer} {Programming} {Problems} with {Umpire}.
\newblock \emph{Management Science}, 20\penalty0 (5):\penalty0 736--773, January 1974.
\newblock ISSN 0025-1909.
\newblock \doi{10.1287/mnsc.20.5.736}.
\newblock URL \url{https://pubsonline.informs.org/doi/10.1287/mnsc.20.5.736}.

\bibitem[Gilpin and Sandholm(2011)]{GILPIN2011147}
Andrew Gilpin and Tuomas Sandholm.
\newblock Information-theoretic approaches to branching in search.
\newblock \emph{Discrete Optimization}, 8\penalty0 (2):\penalty0 147--159, 2011.
\newblock ISSN 1572-5286.
\newblock \doi{https://doi.org/10.1016/j.disopt.2010.07.001}.
\newblock URL \url{https://www.sciencedirect.com/science/article/pii/S1572528610000423}.

\bibitem[Grossmann and Harjunkoski(2019)]{GROSSMANN2019474}
Ignacio~E. Grossmann and Iiro Harjunkoski.
\newblock Process systems engineering: Academic and industrial perspectives.
\newblock \emph{Computers \& Chemical Engineering}, 126:\penalty0 474--484, 2019.
\newblock ISSN 0098-1354.
\newblock \doi{https://doi.org/10.1016/j.compchemeng.2019.04.028}.
\newblock URL \url{https://www.sciencedirect.com/science/article/pii/S009813541930078X}.

\bibitem[Gurobi~Optimization(2020)]{gurobi}
LLC Gurobi~Optimization.
\newblock Gurobi optimizer reference manual, 2020.
\newblock URL \url{http://www.gurobi.com}.

\bibitem[Jezowski(2010)]{jezowski2010review}
Jacek Jezowski.
\newblock Review of water network design methods with literature annotations.
\newblock \emph{Industrial \& Engineering Chemistry Research}, 49\penalty0 (10):\penalty0 4475--4516, 2010.

\bibitem[Kallrath(2021)]{kallrath_optimization:_2021}
Josef Kallrath.
\newblock Optimization: {Using} {Models}, {Validating} {Models}, {Solutions}, {Answers}.
\newblock In Josef Kallrath, editor, \emph{Business {Optimization} {Using} {Mathematical} {Programming}: {An} {Introduction} with {Case} {Studies} and {Solutions} in {Various} {Algebraic} {Modeling} {Languages}}, International {Series} in {Operations} {Research} \& {Management} {Science}, pages 1--32. Springer International Publishing, Cham, 2021.
\newblock ISBN 9783030732370.
\newblock \doi{10.1007/978-3-030-73237-0_1}.
\newblock URL \url{https://doi.org/10.1007/978-3-030-73237-0_1}.

\bibitem[Kim and Dandurand(2022)]{kim2022scalable}
Kibaek Kim and Brian Dandurand.
\newblock Scalable branching on dual decomposition of stochastic mixed-integer programming problems.
\newblock \emph{Mathematical Programming Computation}, 14\penalty0 (1):\penalty0 1--41, 2022.

\bibitem[Kiwiel(1984)]{kc1984algorithm}
Krzysztof~C Kiwiel.
\newblock An algorithm for linearly constrained convex nondifferentiable minimization problems.
\newblock \emph{J. Math. Anal. Appl.(to appear)}, 1984.

\bibitem[Kiwiel(1985)]{kiwiel1985exact}
Krzysztof~C Kiwiel.
\newblock An exact penalty function algorithm for non-smooth convex constrained minimization problems.
\newblock \emph{IMA Journal of Numerical Analysis}, 5\penalty0 (1):\penalty0 111--119, 1985.

\bibitem[Kiwiel(1987)]{kiwiel1987constraint}
Krzysztof~C Kiwiel.
\newblock A constraint linearization method for nondifferentiable convex minimization.
\newblock \emph{Numerische Mathematik}, 51:\penalty0 395--414, 1987.

\bibitem[Kiwiel(1990)]{kiwiel1990proximity}
Krzysztof~C Kiwiel.
\newblock Proximity control in bundle methods for convex nondifferentiable minimization.
\newblock \emph{Mathematical programming}, 46\penalty0 (1-3):\penalty0 105--122, 1990.

\bibitem[Kocis and Grossmann(1987)]{kocis1987relaxation}
Gary~R Kocis and Ignacio~E Grossmann.
\newblock Relaxation strategy for the structural optimization of process flow sheets.
\newblock \emph{Industrial \& engineering chemistry research}, 26\penalty0 (9):\penalty0 1869--1880, 1987.

\bibitem[K{\"u}{\c{c}}{\"u}kyavuz and Sen(2017)]{kuccukyavuz2017introduction}
Simge K{\"u}{\c{c}}{\"u}kyavuz and Suvrajeet Sen.
\newblock An introduction to two-stage stochastic mixed-integer programming.
\newblock In \emph{Leading Developments from INFORMS Communities}, pages 1--27. INFORMS, 2017.

\bibitem[Lee and Leyffer(2011)]{lee_mixed_2011}
Jon Lee and Sven Leyffer.
\newblock \emph{Mixed {Integer} {Nonlinear} {Programming}}.
\newblock Springer Science \& Business Media, December 2011.
\newblock ISBN 9781461419273.

\bibitem[Lemar{\'e}chal(1974)]{lemarechal1974algorithm}
Claude Lemar{\'e}chal.
\newblock An algorithm for minimizing convex functions.
\newblock In \emph{IFIP Congress}, pages 552--556, 1974.

\bibitem[Liao et~al.(2019)Liao, Yao, Han, Fang, Ai, Wen, and He]{LIAO2019265}
Shiwu Liao, Wei Yao, Xingning Han, Jiakun Fang, Xiaomeng Ai, Jinyu Wen, and Haibo He.
\newblock An improved two-stage optimization for network and load recovery during power system restoration.
\newblock \emph{Applied Energy}, 249:\penalty0 265--275, 2019.
\newblock ISSN 0306-2619.
\newblock \doi{https://doi.org/10.1016/j.apenergy.2019.04.176}.
\newblock URL \url{https://www.sciencedirect.com/science/article/pii/S0306261919308359}.

\bibitem[Linderoth and Savelsbergh(1999)]{linderoth1999computational}
Jeff~T Linderoth and Martin~WP Savelsbergh.
\newblock A computational study of search strategies for mixed integer programming.
\newblock \emph{INFORMS Journal on Computing}, 11\penalty0 (2):\penalty0 173--187, 1999.

\bibitem[McCormick(1976)]{mccormic_1976}
Garth~P. McCormick.
\newblock Computability of global solutions to factorable nonconvex programs: Part i ---convex underestimating problems.
\newblock \emph{Mathematical Programming}, 10\penalty0 (1):\penalty0 147--175, 1976.
\newblock \doi{10.1007/BF01580665}.
\newblock URL \url{https://doi.org/10.1007/BF01580665}.

\bibitem[Misener and Floudas(2012)]{misener2012global}
Ruth Misener and Christodoulos~A Floudas.
\newblock Global optimization of mixed-integer quadratically-constrained quadratic programs (miqcqp) through piecewise-linear and edge-concave relaxations.
\newblock \emph{Mathematical Programming}, 136\penalty0 (1):\penalty0 155--182, 2012.

\bibitem[Morrison et~al.(2016)Morrison, Jacobson, Sauppe, and Sewell]{MORRISON201679}
David~R. Morrison, Sheldon~H. Jacobson, Jason~J. Sauppe, and Edward~C. Sewell.
\newblock Branch-and-bound algorithms: A survey of recent advances in searching, branching, and pruning.
\newblock \emph{Discrete Optimization}, 19:\penalty0 79--102, 2016.
\newblock ISSN 1572-5286.
\newblock \doi{https://doi.org/10.1016/j.disopt.2016.01.005}.
\newblock URL \url{https://www.sciencedirect.com/science/article/pii/S1572528616000062}.

\bibitem[Murtagh and Saunders(1995)]{murtagh1995minos}
Bruce~A Murtagh and Michael~A Saunders.
\newblock Minos 5.4 user's guide (revised).
\newblock Technical report, Technical Report SOL 83-20R, Department of Operations Research, Stanford~…, 1995.

\bibitem[Oliveira and Sagastiz{\'a}bal(2014)]{oliveira2014bundle}
Welington~de Oliveira and Claudia Sagastiz{\'a}bal.
\newblock Bundle methods in the xxist century: A bird's-eye view.
\newblock \emph{Pesquisa Operacional}, 34:\penalty0 647--670, 2014.

\bibitem[Ortega and Wolsey(2003)]{ortega2003branch}
Francisco Ortega and Laurence~A Wolsey.
\newblock A branch-and-cut algorithm for the single-commodity, uncapacitated, fixed-charge network flow problem.
\newblock \emph{Networks: An International Journal}, 41\penalty0 (3):\penalty0 143--158, 2003.

\bibitem[Palani et~al.(2019)Palani, Wu, and Morcos]{palani2019frank}
Ananth~M Palani, Hongyu Wu, and Medhat~M Morcos.
\newblock A frank--wolfe progressive hedging algorithm for improved lower bounds in stochastic scuc.
\newblock \emph{IEEE Access}, 7:\penalty0 99398--99406, 2019.

\bibitem[Rockafellar and Wets(1991)]{rockafellar1991scenarios}
R~Tyrrell Rockafellar and Roger J-B Wets.
\newblock Scenarios and policy aggregation in optimization under uncertainty.
\newblock \emph{Mathematics of operations research}, 16\penalty0 (1):\penalty0 119--147, 1991.

\bibitem[Sewell et~al.(2012)Sewell, Sauppe, Morrison, Jacobson, and Kao]{sewell_bbr_2012}
Edward~C. Sewell, Jason~J. Sauppe, David~R. Morrison, Sheldon~H. Jacobson, and Gio~K. Kao.
\newblock A bb\&r algorithm for minimizing total tardiness on a single machine with sequence dependent setup times.
\newblock \emph{Journal of Global Optimization}, 54\penalty0 (4):\penalty0 791--812, 2012.
\newblock \doi{10.1007/s10898-011-9793-z}.
\newblock URL \url{https://doi.org/10.1007/s10898-011-9793-z}.

\bibitem[The Optimization~Firm(2019)]{baron}
LLC The Optimization~Firm.
\newblock Baron user manual v. 2019.12.7, 2019.
\newblock URL \url{https://www.minlp.com}.

\bibitem[Vilà and Pereira(2014)]{VILA2014105}
Mariona Vilà and Jordi Pereira.
\newblock A branch-and-bound algorithm for assembly line worker assignment and balancing problems.
\newblock \emph{Computers \& Operations Research}, 44:\penalty0 105--114, 2014.
\newblock ISSN 0305-0548.
\newblock \doi{https://doi.org/10.1016/j.cor.2013.10.016}.
\newblock URL \url{https://www.sciencedirect.com/science/article/pii/S0305054813003110}.

\bibitem[Virasjoki et~al.(2020)Virasjoki, Siddiqui, Oliveira, and Salo]{virasjoki2020utility}
Vilma Virasjoki, Afzal~S Siddiqui, Fabricio Oliveira, and Ahti Salo.
\newblock Utility-scale energy storage in an imperfectly competitive power sector.
\newblock \emph{Energy Economics}, 88:\penalty0 104716, 2020.

\bibitem[Wang et~al.(2021)Wang, Meng, Wang, and Qu]{WANG2021143}
Tingsong Wang, Qiang Meng, Shuaian Wang, and Xiaobo Qu.
\newblock A two-stage stochastic nonlinear integer-programming model for slot allocation of a liner container shipping service.
\newblock \emph{Transportation Research Part B: Methodological}, 150:\penalty0 143--160, 2021.
\newblock ISSN 0191-2615.
\newblock \doi{https://doi.org/10.1016/j.trb.2021.04.016}.
\newblock URL \url{https://www.sciencedirect.com/science/article/pii/S0191261521000795}.

\bibitem[Watson and Woodruff(2011)]{watson2011progressive}
Jean-Paul Watson and David~L Woodruff.
\newblock Progressive hedging innovations for a class of stochastic mixed-integer resource allocation problems.
\newblock \emph{Computational Management Science}, 8\penalty0 (4):\penalty0 355--370, 2011.

\bibitem[wen Chen(2003)]{CHEN20031925}
Xue wen Chen.
\newblock An improved branch and bound algorithm for feature selection.
\newblock \emph{Pattern Recognition Letters}, 24\penalty0 (12):\penalty0 1925--1933, 2003.
\newblock ISSN 0167-8655.
\newblock \doi{https://doi.org/10.1016/S0167-8655(03)00020-5}.
\newblock URL \url{https://www.sciencedirect.com/science/article/pii/S0167865503000205}.

\bibitem[Zhao and Luh(2002)]{zhao2002new}
X.~Zhao and P.B. Luh.
\newblock New bundle methods for solving lagrangian relaxation dual problems.
\newblock \emph{Journal of Optimization Theory and Applications}, 113\penalty0 (2):\penalty0 373--397, May 2002.
\newblock ISSN 0022-3239, 1573-2878.
\newblock \doi{10.1023/A:1014839227049}.
\newblock URL \url{http://link.springer.com/10.1023/A:1014839227049}.

\end{thebibliography}

\appendix
\section{Appendix A}
\label{sec:ap1}

\begin{table}[!h]
\caption{Range of the parameters' values for RDEM problem. The values utilised are uniformly sampled from these ranges.}
\centering
\begin{tabular}{cc}
\toprule
Parameter & Range      \\
\hline
$Q$ & [-100, 100]    \\
$c$ & [0, 100]    \\
$q$ & [-100, 100]    \\
$U$ & [0, 1000]   \\
$T$ & [0, 100]    \\
$W$ & [-100, 100]    \\
$h$ & [1000, 100000] \\
$X$ & [0, 10]     \\
$Y$ & [0, 10] \\ 
\bottomrule
\end{tabular}
\end{table}

\section{Appendix B}
\label{sec:ap2}

\begin{table}[!h]
\caption{Average dual bound generated
within one hour by solving the Instances 1-5 with Gurobi solver}
\centering
\begin{tabular}{cc}
\toprule
Instance &  dual bound (\%)     \\
\hline
1 & 0.00    \\ 
2 & 4.71  \\ 
3 & 6.52   \\ 
4 & 8.54   \\ 
5 & 5.11   \\ 
\bottomrule
\end{tabular}
\end{table}

\section{Appendix C}
\label{sec: ap3}

This section presents an illustrative example showcasing the application of Algorithm \ref{al: dual_bnb}. Consider Problem \eqref{eq: RDEM} with 2 scenarios $S = \{s_1, s_2\}$, two first-stage variables $x_1$ and $x_2$, one second-stage variable $y$, and one constraint defined per scenario. Then,  Problem \eqref{eq: RDEM} can be formulated as follows: 

\begin{align}
    z^\text{SMIP} = \mini_{x_{1}^{1}, x_{1}^{2}, x_{2}^{1}, x_{2}^{2}, y^{1}, y^{2}} 
    ~&  -16.8 y^{1}y^{1} - 9.8 y^{1}  \\  ~& -4.5 y^{2}y^{2} - 7.9 y^{2} \\  
    ~& -53.2 x_1^1  -23.2 x_2^1 - 9.0 x_1^2 - 3.9 x_2^2 \\ 
    \text{s.t.:} \nonumber \\ 
    ~&  78.1 x_1^1 + 67.0 x_2^1 + 16.8 y^1 \le 57108.7 \\ 
    ~& 77.0 y^2y^2 + 45.3 x_1^2 + 30.2 x_2^2 + 0.1 y^2  \le 56702.4 \\
    ~& x_1^1 - \hat{x}_1 = 0 \label{il: na_1} \\ 
    ~& x_1^2 - \hat{x}_1 = 0 \label{il: na_2} \\ 
    ~& x_2^1 - \hat{x}_2 = 0  \label{il: na_3} \\
    ~& x_2^2 - \hat{x}_2 = 0  \label{il: na_4} \\
     ~& 0.0 \le x_1^1, x_1^2 \le 5.0 \\
    ~& 0.0 \le  x_2^1, x_2^2 \le  7.0 \\
 ~& 0.0 \le y^1 \le 9.0 \\
  ~& 0.0 \le y^2 \le 8.0 \\
     ~& x_1^1, \ x_1^2, \ x_2^1, \ x_2^2 \  \in \  \mathbb{Z} \\
\end{align} 

Where, \eqref{il: na_1}- \eqref{il: na_4} are non-anticipativity conditions.

Applying \textit{p}-Lagrangian relaxation we get Problem \eqref{p_LR} with two sub-problems \eqref{il: plr_1} and \eqref{il: plr_2}

\begin{equation}
\label{il: plr_1}
\begin{aligned}
   L^1(\mu^1) = \mini_{\omega^1} 
    ~& \frac{1}{\pi^1} \left(-16.8 w^1 - 53.3 x_1^1 - 23.2 x_2^1 - 9.8 y^1 + \mu^1_1 (x_1^1-\hat{x}_1) + \mu^1_2 (x_2^1-\hat{x}_2) \right)\\ 
    \text{s.t.:}  \\
    ~& y^1 - 9 \Delta y^1 - 4.5 z_{-1}^1 = 0 \\ 
    ~& -4.5 \hat{y}_{-1}^1 + w^1 - 9 \Delta w^1 = 0 \\ 
    ~& 78.1 x_1^1 + 67.0 x_2^1 + 16.8 y^1 \leq 57108.7 \\
    ~& 0.5 y^1 + 9 \Delta y^1 - \Delta w^1 \leq 4.5 \\
    ~& -0.5 y^1 + \Delta w^1 \leq 0 \\
    ~& -\Delta w^1 \leq 0 \\
    ~& -9 \Delta y^1 + \Delta w^1 \leq 0 \\
    ~& -\hat{y}_{-1}^1 \leq 0 \\
    ~& \hat{y}_{-1}^1 - 9 z_{-1}^1 \leq 0 \\
    ~& -y^1 + \hat{y}_{-1}^1 \leq 0 \\
    ~& y^1 - \hat{y}_{-1}^1 + 9 z_{-1}^1 \leq 9.0 \\
    ~& 0.0 \leq x_1^1 \leq 5.0 \\
    ~& 0.0 \leq x_2^1 \leq 7.0 \\
    ~& 0.0 \leq y^1 \leq 9.0 \\
    ~& 0.0 \leq \Delta y^1 \leq 0.5 \\
    ~& x_1^1, x_2^1 \in \mathbb{Z} \\
    ~& z_{-1}^1 \in \{0,1\}
\end{aligned}
\end{equation}

\begin{equation}
\label{il: plr_2}
\begin{aligned}
    L^2(\mu^2) = \mini_{\omega^2} 
    ~& \frac{1}{\pi^2} \left(-4.5 w^2 - 9.0 x_1^2 - 3.9 x_2^2 - 7.9 y^2 + \mu^2_1 (x_1^2-\hat{x}_1) + \mu^2_1 (x_2^2-\hat{x}_2) \right)\\ 
    \text{s.t.:}\\
    ~& y^2 - 8 \Delta y^2 - 4 z_{-1}^2 = 0 \\ 
    ~& -4 \hat{y}_{-1}^2 + w^2 - 8 \Delta w^2 = 0 \\ 
    ~& 45.3 x_1^2 + 30.2 x_2^2 + 0.1 y^2 + 77 w^2 \leq 56702.4 \\
    ~& 0.5 y^2 + 8 \Delta y^2 - \Delta w^2 \leq 4.0 \\
    ~& -0.5 y^2 + \Delta w^2 \leq 0 \\
    ~& -\Delta w^2 \leq 0 \\
    ~& -8 \Delta y^2 + \Delta w^2 \leq 0 \\
    ~& -\hat{y}_{-1}^2 \leq 0 \\
    ~& \hat{y}_{-1}^2 - 8 z_{-1}^2 \leq 0 \\
    ~& -y^2 + \hat{y}_{-1}^2 \leq 0 \\
    ~& y^2 - \hat{y}_{-1}^2 + 8 z_{-1}^2 \leq 8.0 \\
    ~& 0.0 \leq x_1^2 \leq 5.0 \\
    ~& 0.0 \leq x_2^2 \leq 7.0 \\
    ~& 0.0 \leq y^2 \leq 8.0 \\
    ~& 0.0 \leq \Delta y^2 \leq 0.5 \\
    ~& x_1^2, x_2^2 \in \mathbb{Z} \\
    ~& z_{-1}^2 \in \{0,1\}
\end{aligned}
\end{equation}

The terms $\pi^1$ and $\pi^2$ correspond to the probabilities of Problems \eqref{il: plr_1} and \eqref{il: plr_2}, respectively, and  $$\omega^1 = \{x_1^1, x_2^1, \hat{x_1}, \hat{x_2}, y^1, \Delta y^1, \Delta w^1, w^1, \hat{y}_{-1}^1, z_{-1}^1\}, $$ and
$$\omega^2 = \{x_1^2, x_2^2, \hat{x_1}, \hat{x_2}, y^2, \Delta y^2, \Delta w^2, w^2, \hat{y}_{-1}^2, z_{-1}^2\}.$$

Then the first root node of the Algorithm \ref{al: dual_bnb} will be 

\begin{align}
\label{il: root node}
z_{LD} = \maxi_{\mu^1, \mu^2} \left\{ \pi^1 L^1(\mu^1) + \pi^2 L^2(\mu^2): \pi^1 \mu^1 + \pi^2 \mu^2 = 0\right\}.
\end{align}

Applying Algorithm \ref{al: proximal bundle method} to solve Problem \eqref{il: root node} yields the solution 
in just one iteration, due to the simple structure and small scale of the problem. A similar result occurs when applying Algorithm \ref{al: FWPH}, which also converges to the solution 
$z_{LB}$ in one iteration. The solution is as follows:

\begin{table}[h]
\label{il: bm}
\caption{Result of applying Algorithm \ref{al: proximal bundle method} and Algorithm \ref{al: FWPH} to Problem \eqref{il: root node} }
\centering
\begin{tabular}{ll}
\hline
Number of iterations   &  1 \\
Objective function value      & $z_{LB}= -2299.62$ \\
First-stage variables' values  & $x_1^1 = 5.0$, $x_2^1 = 7.0$, $x_1^2 = 5.0$, $x_2^2 = 7.0$ \\
Second-stage variables' values & $y^1 = 9.0$, $y^2 = 8.0$ \\ \hline
\end{tabular}
\end{table}

The average solution is given by $\hat{x} = \pi^1 x^1 + \pi^2 x^2$, resulting in $\{5, 7\}$. Due to the simplicity of the problem, both the \eqref{al: proximal bundle method} and Algorithm \ref{al: FWPH} managed to yield a solution at the root node that satisfies both the non-anticipativity and integrality conditions ($x_1^1 = x_1^2 = 5.0$ and $x_2^1 = x_2^2 = 7.0$). However, in the general case, there is no guarantee that this will occur. Therefore, for the sake of illustration, let us assume that the solution to Problem \eqref{il: root node} violates one of the conditions, considering $\pi^1 = \pi^2 = 0.5$.

First, let us assume that $x_1^1 = 4$, $x_1^2 = 5$, and $x_2^1 = x_2^2 = 7$. In this case, the average solution is $\hat{x} = \{4.5, 7\}$. We can observe that the first component takes the value 4.5, which is not integer-valued. In this situation, following Algorithm \ref{al: dual_bnb}, we would perform branching based on the violation of the integrality condition. Specifically, we would generate two child nodes $N_1$ and $N_2$, such that $N_1$ contains Problem \eqref{il: root node} with the additional constraints $\{x_1^1, x_1^2 \le 4 \}$, while $N_2$ contains Problem \eqref{il: root node} with the constraints $\{x_1^1, x_1^2 \ge 5 \}$.

Alternatively, let us assume that $x_1^1 = 4$, $x_1^2 = 6$, and $x_2^1 = x_2^2 = 7$. In this case, the average solution is $\hat{x} = \{5, 7\}$. We can see that the solution $\hat{x}$ satisfies the integrality condition but does not satisfy the non-anticipativity condition for the first component. Specifically, the dispersion is $\sigma_1 = x_1^2 - x_1^1 = 2$. In this case, following Algorithm \ref{al: dual_bnb}, we would perform branching based on the violation of the non-anticipativity condition. We would create two child nodes $N_1$ and $N_2$, where the subproblem for $N_1$ corresponds to Problem \eqref{il: root node} with the constraints $\{x_1^1, x_1^2 \le 5 - \epsilon_{BB}\}$, and $N_2$ would correspond to Problem \eqref{il: root node} with the constraints $\{x_1^1, x_1^2  \ge 5 + \epsilon_{BB}\}$.

\end{document}